\documentclass[a4paper,11pt]{article}

\usepackage{amsmath,amssymb}
\usepackage[T1]{fontenc}
\usepackage{latexsym}
\usepackage[amsmath,thmmarks]{ntheorem}
\usepackage{enumerate}
\usepackage{fancyhdr}
\usepackage{theorem}
\usepackage{scrpage}
\usepackage{ucs}
\usepackage{latexsym}
\usepackage{dsfont}				
\usepackage{mathrsfs}			
\usepackage[active]{srcltx}		
\usepackage{enumitem}

\title{Paraproducts via $H^{\infty}$-functional calculus}
\author{Dorothee Frey}

\DeclareMathOperator{\supp}{supp}
\DeclareMathOperator{\dist}{dist}

\DeclareMathOperator{\loc}{loc}
\numberwithin{equation}{section}

\newtheorem{Def}{Definition}[section]
\newtheorem{Lemma}[Def]{Lemma}

\newtheorem{Theorem}[Def]{Theorem}
\newtheorem{Prop}[Def]{Proposition}
{\theorembodyfont{\rmfamily} 
}
{\theorembodyfont{\rmfamily} \newtheorem{Remark}[Def]{Remark}}
{\theoremstyle{nonumberplain}
\theorembodyfont{\upshape}
\theoremheaderfont{\normalfont \bfseries}
\theoremsymbol{\ensuremath{_\square}}
\theoremseparator{:}
\newtheorem{Proof}{Proof}}

\newcommand{\R}{\ensuremath{\mathbb{R}}}
\newcommand{\C}{\ensuremath{\mathbb{C}}}
\newcommand{\N}{\ensuremath{\mathbb{N}}}

\newcommand{\bH}{\ensuremath{\mathbb{H}}}
\newcommand{\Eins}{\ensuremath{\mathds{1}}}

\newcommand{\calD}{\ensuremath{\mathcal{D}}}

\newcommand{\calE}{\ensuremath{\mathcal{E}}}
\newcommand{\calM}{\ensuremath{\mathcal{M}}}
\newcommand{\calN}{\ensuremath{\mathcal{N}}}
\newcommand{\calC}{\ensuremath{\mathcal{C}}}
\newcommand{\calR}{\ensuremath{\mathcal{R}}}

\newcommand{\calO}{\ensuremath{\mathcal{O}}}
\newcommand{\scrC}{\ensuremath{\mathscr{C}}}
\newcommand{\scrA}{\ensuremath{\mathscr{A}}}
\newcommand{\norm}[1]{\left\|#1\right\|}
\newcommand{\abs}[1]{\left|#1\right|}
\newcommand{\abbs}[1]{|#1|}
\newcommand{\nnorm}[1]{\|#1\|}
\newcommand{\skp}[1]{\langle #1 \rangle}

\newcommand{\eps}{\ensuremath{\varepsilon}}

\newcommand{\molMeps}{\calM_0^{1,2,M,\eps}(L)}

\newcommand{\molMdeps}{\calM_0^{1,2,M,\eps}(L^{\ast})}

\begin{document}

\maketitle

\begin{abstract}
Let $X$ be a space of homogeneous type and let $L$ be a sectorial operator with bounded holomorphic functional calculus on $L^2(X)$. We assume that the semigroup $\{e^{-tL}\}_{t>0}$ satisfies Davies-Gaffney estimates. 
In this paper, we introduce a new type of paraproduct operators that is constructed via certain approximations of the identity associated to $L$. 
We show various boundedness properties on $L^p(X)$ and the recently developed Hardy and BMO spaces $H^p_L(X)$ and $BMO_L(X)$.
In generalization of standard paraproducts constructed via convolution operators, we show  $L^2(X)$ off-diagonal estimates as a substitute for Calder\'on-Zygmund kernel estimates.
As an application, we study differentiability properties of paraproducts in terms of fractional powers of the operator $L$. \\
The results of this paper are fundamental for the proof of a $T(1)$-Theorem for operators beyond Calder\'on-Zygmund theory, which will be the subject of a forthcoming paper.
\end{abstract}

\tableofcontents

\section{Introduction and main results}

Paraproduct operators are an important tool in harmonic analysis and play an essential role in analysis and the theory of partial differential equations. They emerged in the theory of paradifferential operators, see e.g. \cite{CoifmanMeyer2} and \cite{Bony}, and have crucial applications in the general theory of singular integral operators and the study of non-linear problems, see e.g. \cite{KatoPonce} in the context of Euler and Navier-Stokes equations.\\ 
More specific, in the proof of the $T(1)$-Theorem of David and Journ\'e \cite{DavidJourne}, a main ingredient is the following paraproduct. 
Given $b \in BMO(\R^n)$, one defines an operator $\Pi_b$  
on $L^2(\R^n)$ via
\begin{align} \label{Para-DJ}
 	\Pi_b f = \int_0^{\infty} Q_t[(Q_t b) (P_t f)] \,\frac{dt}{t}, \qquad f \in L^2(\R^n),
\end{align}
where $P_t$ and $Q_t$  are convolution operators with $P_t(1)=1$ and $Q_t(1)=0$.
One can then show that $\Pi_b$ is a Calder\'{o}n-Zygmund operator, bounded on $L^2(\R^n)$ and satisfying $\Pi_b(1)=b$ and $\Pi_b^{\ast}(1)=0$. \\

In the last two decades, the study of properties of sectorial operators often depended on pointwise Gaussian estimates for the kernel of the corresponding semigroup, which therefore acts bounded on $L^p$ for $p \in [1,\infty]$. 
In recent years, there was then developed a theory for sectorial operators $L$ whose semigroup is bounded on $L^p$ only for a range of $p$ strictly smaller than $(1,\infty)$.
For such operators, one cannot work with pointwise Gaussian estimates for the semigroup, but has to work with generalized Gaussian estimates, Davies-Gaffney estimates or other off-diagonal estimates instead. A key role in this theory is played by approximation operators that are constructed via $H^{\infty}$-functional calculus as introduced in \cite{McIntosh}. For example, the semigroup $\{e^{-tL}\}_{t>0}$ can be used as an approximation of the identity and the derivative $\{t\partial_te^{-tL}\}_{t>0}$ for the construction of a resolution identity. 
In this way, there were obtained various results on generalizations of operators and function spaces, that were originally constructed via the Laplacian and Littlewood-Paley theory. This includes Hardy spaces $H^p_L$ and a corresponding space $BMO_L$ that are associated to $L$, see e.g. \cite{AuscherDuongMcIntosh}, \cite{DuongYan2}, \cite{AuscherMcIntoshRuss}, \cite{BernicotZhao}, \cite{HofmannMayboroda}, \cite{HofmannMayborodaMcIntosh}, \cite{HLMMY}, \cite{DuongLi}, Riesz transforms, e.g. in \cite{KatoSquare}, \cite{HofmannMartell}, \cite{BlunckKunstmann2}, and other operators beyond Calder\'on-Zygmund theory, e.g. in \cite{BlunckKunstmann}, \cite{ACDH}, \cite{AuscherD}, \cite{Auscher}.\\

In this article, we introduce a new type of paraproduct operator and generalize the above paraproduct in the following sense.\\
We assume $X$ to be a space of homogeneous type and let $L$ be a sectorial operator with bounded holomorphic functional calculus on $L^2(X)$. We assume that the semigroup $\{e^{-tL}\}_{t>0}$ satisfies Davies-Gaffney estimates and, for some results, an $L^p-L^2$ estimate for some $p<2$. 
Standard examples of operators that satisfy our assumptions are elliptic operators in divergence form with bounded complex coefficients, see e.g. \cite{Auscher}, Schr\"odinger operators with singular potentials, see e.g. \cite{LSV}, and Laplace-Beltrami operators on complete Riemannian manifolds with non-negative Ricci curvature, see e.g. \cite{Davies1}, \cite{Grigoryan}.\\
With help of the $H^{\infty}$-functional calculus, we define a paraproduct associated to $L$ by
 \begin{equation} \label{intro-para}
 	\Pi_b: f \mapsto  \int_0^{\infty} \tilde{\psi}(t^{2m}L) [\psi(t^{2m}L)b \cdot A_t(e^{-t^{2m}L}f)] \, \frac{dt}{t},
 \end{equation}
where $\psi$, $\tilde{\psi}$ are taken from the set $\Psi$ consisting of bounded holomorphic functions on a sector with decay at zero and infinity, e.g. $\psi(tL)=(tL)^Me^{-tL}$ for $M>\frac{n}{4m}$, and $A_t$ denotes some averaging operator. \\
The appearance of the operator $A_t$ might seem to be surprising, but this is due to the fact that we do not impose any kernel estimates on the semigroup $\{e^{-tL}\}_{t>0}$.\\
For $X=\R^n$ and $L=-\Delta$, one can omit the averaging operator $A_t$ and the definition in \eqref{intro-para} then corresponds to paraproducts defined via convolution.\\

Paraproducts defined in this way allow for a great flexibility, making it possible to adapt them to many situations in Calder\'on-Zygmund theory, and, more importantly, beyond Calder\'on-Zygmund theory.
The spaces $H^p_L(X)$ and $BMO_L(X)$, that are associated to $L$, generalize the usual Lebesgue spaces and the space BMO of John and Nirenberg and are the appropriate setting for paraproducts of the form \eqref{intro-para}.\\
Our first main result is the following. 

\begin{Theorem} \label{Thm1-intro}
 Let $b \in BMO_L(X)$ and let $\psi,\tilde{\psi}$ as specified in Theorem \ref{paraproduct}. Then $\Pi_b$, defined in \eqref{intro-para}, is bounded on $L^2(X)$ and extends to a bounded operator from $L^p(X)$ to $H^p_L(X)$ for $p \in (2,\infty)$ and from $L^{\infty}(X)$ to $BMO_L(X)$.
\end{Theorem}

Moreover, the conservation property $e^{-tL}(1)=1$ in $L^2_{\loc}(X)$ allows us to reobtain the properties $\Pi_b(1)=b$ and $\Pi_b^{\ast}(1)=0$. \\
For a second order elliptic operator $L$ in divergence form, we denote by $(p_{-}(L),p_{+}(L))$ the interior of the interval of $L^p$ boundedness of $\{e^{-tL}\}_{t>0}$. Then for $p \in  (p_{-}(L),p_{+}(L))$, as shown in \cite{HofmannMayborodaMcIntosh}, there holds $H^p_L(X)=L^p(X)$, and therefore $\Pi_b$ is bounded on $L^p(X)$ for all $p \in [2,p_{+}(L))$. For other types of operators $L$, one can obtain similar results via generalized Gaussian estimates, cf. Proposition \ref{Hp-equiv} below. \\
The proof of Theorem \ref{paraproduct} heavily relies on an analogue of the Fefferman-Stein criterion.
That is, except for a growth estimate for $b$, there holds
\[
		b \in BMO_L(X) \quad \Longleftrightarrow \quad \nu_{\psi,b} := \abs{\psi(t^{2m}L)b(y)}^2 \frac{d\mu(y)\, dt}{t} \quad \text{is a Carleson measure}.
\]
For $\psi(z)=z^Me^{-z}$, $M>\frac{n}{4m}$, the result is due to \cite{HofmannMayboroda}. We generalize the result to allow for a greater freedom in the choice of $\psi$, cf.  Proposition \ref{LemmaCarlesonMeasure}.\\

We continue our studies of paraproducts by defining $\Pi(f,b):=\Pi_b(f)$ and considering the paraproduct as an operator of the second variable. 
In analogy to the fact that the paraproduct in \eqref{Para-DJ} is a Calder\'{o}n-Zygmund operator, we show certain off-diagonal estimates for the paraproduct associated to $L$. These off-diagonal estimates, that have been used before,  e.g. in \cite{BlunckKunstmann}, \cite{Auscher}, \cite{HofmannMayboroda}, enable us to extend the operator to certain $L^p(X)$ and $H^p_L(X)$ spaces. We obtain the following result.

\begin{Theorem} \label{Thm2-intro}
Let  $\psi,\tilde{\psi}$ as specified in Theorem \ref{para-bdd-hp}.
Then $\Pi: L^{\infty}(X) \times L^2(X) \to L^2(X)$ is bounded and extends to a bounded operator
$\Pi: L^{\infty}(X) \times H^p_L(X) \to L^p(X)$ for  $p \in [1,2)$ and 
$\Pi: L^{\infty}(X) \times L^p(X) \to H^p_L(X)$ for $p \in (2,\infty)$.
\end{Theorem}

As before, the identification of $H^p_L(X)$ and $L^p(X)$ for a certain range of $p$ (according to \cite{HofmannMayborodaMcIntosh} or to Proposition \ref{Hp-equiv} below) yields boundedness results $\Pi: L^{\infty}(X) \times L^p(X) \to L^p(X)$.\\ 

We end the article with some results on differentiability properties of paraproducts constructed via $H^{\infty}$-functional calculus and show that there holds a Leibniz-type rule. More results will be given in \cite{FreyHytoenen}.\\
An important application of the paraproduct defined in \eqref{intro-para} is given in \cite{FreyKunstmann}, where we generalize the $T(1)$-Theorem for operators beyond Calder\'on-Zygmund theory.\\
While the work was in preparation, we learned that similar paraproducts have also been considered by Bernicot, cf. \cite{Bernicot}. The main difference to our results is, that a crucial assumption in \cite{Bernicot} are pointwise bounds on the kernels of the semigroup $\{e^{-tL}\}_{t>0}$, an assumption which is considerably relaxed here.\\

The article is organized as follows: In Section 2 we collect the most important definitions and results of $H^{\infty}$-functional calculus, tent spaces and Carleson measures and fix our assumptions on the operator $L$. In Section 3 we unify the theory of Hardy and BMO spaces associated to operators. We generalize the results, in the literature so far only stated for second order operators, to higher order operators and prove a generalization of a Calder\'on reproducing formula and a Carleson measure characterization of $BMO_L(X)$.
Section 4 is devoted to statement and proof of our main results, Theorem \ref{Thm1-intro} and Theorem \ref{Thm2-intro}. We end with a Leibniz-type rule.\\

Throughout the article, the letter ``$C$'' will denote (possibly different) positive constants that are independent of the essential variables. We will frequently write $a \lesssim b$ for non-negative quantities $a,b$, if $a \leq Cb$ for some $C$.


\section{Preliminaries}

In the following we will always assume $X$ to be a space of homogeneous type. 
More precisely, we assume that $(X,d)$ is a metric space and $\mu$ is a nonnegative Borel measure on $X$ with $\mu(X)=\infty$ which satisfies the \emph{doubling condition}:\\
 There exists a constant $A_1 \geq 1$ such that for all $x \in X$ and all $r>0$
 \begin{align*} 
 		V(x,2r) \leq A_1 V(x,r) < \infty,
 \end{align*}
 where we set $B(x,r):=\{y \in X\,:\, d(x,y)<r\}$ and $V(x,r):=\mu(B(x,r))$.

 Note that the doubling property implies the following strong homogeneity property: There exists a constant $A_2>0$  and some $n>0$ such that for all $\lambda \geq 1$, for all $x \in X$ and all $r>0$
 \begin{align} \label{doublingProperty2}
 		V(x,\lambda r) \leq A_2 \lambda^n V(x,r).
 \end{align}
 In a Euclidean space with the Lebesgue measure, the parameter $n$ corresponds to the dimension of the space.
 For more details on spaces of homogeneous type, see \cite{CoifmanWeiss2}.\\
For a ball $B \subseteq X$ we denote by $r_B$ the radius of $B$ and set
\begin{equation} \label{annuli}
 	S_0(B):=B \qquad \text{and} \qquad S_j(B):= 2^jB \setminus 2^{j-1} B \quad \text{for} \; j=1,2,\ldots,
\end{equation}
where $2^jB$ is the ball with the same center as $B$ and radius $2^jr_B$. \\
Let $t>0$. We define the \emph{averaging operator} $A_t$ by
\begin{align} \label{avOp}
	A_tf(x) := \frac{1}{V(x,t)} \int_{B(x,t)} f(y) \,d\mu(y)
\end{align}
for all $x \in X$ and every $f \in L^1_{\loc}(X)$.\\
We denote by $\calM$ the uncentered Hardy-Littlewood maximal operator. For $p \in [1,\infty)$ and measurable functions $f: X \to \C$ we set $ \calM_p f := [\calM(\abs{f}^p)]^{1/p}$.

\subsection{Holomorphic functional calculus}

We only state the most important definitions and results. For more details on holomorphic functional calculi we refer to \cite{McIntosh}, \cite{AlbrechtDuongMcIntosh}, \cite{LevicoNotes} and \cite{Haase}. \\
For $0 \leq \omega < \sigma < \pi$ we define the closed and open sectors in the complex plane $\C$ by
\begin{align*}
 	S_{\omega+} 	:= \{\zeta \in \C \setminus \{0\}\,:\, \abs{\arg \zeta} \leq \omega\} \cup \{0\}, \qquad 
	\Sigma^0_{\sigma} 		:= \{\zeta \in \C \,:\, \zeta \neq 0, \abs{\arg \zeta} < \sigma\}.
\end{align*}
We denote by $H(\Sigma^0_{\sigma})$ the space of all holomorphic functions on $\Sigma^0_{\sigma}$.
We further define 
\begin{align*}
 	H^{\infty}(\Sigma_{\sigma}^0)		&:= \{\psi \in H(\Sigma^0_{\sigma}) \,:\, \norm{\psi}_{L^{\infty}(\Sigma_{\sigma}^0) }< \infty \}, \\
	\Psi_{\alpha,\beta} (\Sigma^0_{\sigma}) 
				& := \{ \psi \in H(\Sigma^0_{\sigma}) \,:\, \abs{\psi(\zeta)} \leq C \abs{\zeta}^\alpha(1+\abs{\zeta}^{\alpha+\beta})^{-1} \text{ for every } \zeta \in \Sigma_\sigma^0\}
\end{align*}
for every $\alpha,\beta>0$ and $\Psi(\Sigma_{\sigma}^0):= \bigcup_{\alpha,\beta>0} \Psi_{\alpha,\beta}(\Sigma_{\sigma}^0)$.

\begin{Def}
Let $\omega \in [0, \pi)$. A closed operator $L$ in a Hilbert space $H$ is said to be  \emph{sectorial of angle $\omega$}
if $\sigma(L) \subseteq S_{\omega+}$ and, for each $\sigma > \omega$, there exists a constant $C_{\sigma}>0$ such that
\[
 	\norm{(\zeta I-L)^{-1}} \leq C_{\sigma} \abs{\zeta}^{-1}, \qquad \zeta \notin S_{\sigma+}.
\]
\end{Def}

\begin{Remark}
Let  $\omega \in [0, \pi)$ and let $L$ be a sectorial operator of angle $\omega$ in a Hilbert space $H$. Then $L$ has dense domain in $H$. If $L$ is assumed to be injective, then $L$ also has dense range in $H$. See e.g. \cite{CDMcY}, Theorem 2.3 and Theorem 3.8.
\end{Remark}

Let $\omega < \theta < \sigma<\pi$ and let $L$ be a sectorial operator of angle $\omega \in [0,\pi)$ in a Hilbert space $H$. 
Then for every $\psi \in \Psi(\Sigma_{\sigma}^0)$
\begin{align} \label{Def-functcalc}
 	\psi(L):=\frac{1}{2\pi i} \int_{\partial \Sigma_{\theta}^0} \psi(\lambda) (\lambda I-L)^{-1} \,d\lambda
\end{align}
defines a bounded operator on $H$. 
 By sectoriality of $L$ the integral in \eqref{Def-functcalc} is well-defined, and an extension of Cauchy's theorem shows that the definition is independent of the choice of $\theta \in (\omega,\sigma)$.\\
Let $L$ be in addition injective and set $\psi(z):=z(1+z)^{-2}$. Then $\psi(L)$ is injective and has dense range in $H$. For $f \in H^{\infty}(\Sigma_{\sigma}^0)$
 one can define by
 \begin{align*}
 			f(L):= [\psi(L)]^{-1}(f \cdot \psi)(L) 
 \end{align*}
a closed operator in $H$. We say that $L$ has a \emph{bounded} $H^{\infty}(\Sigma_{\sigma}^0)$ \emph{functional calculus} if there exists a constant $c_{\sigma}>0$ such that for all $f \in H^{\infty}(\Sigma^0_{\sigma})$, there holds $f(L) \in B(H)$ with
\[
 	\norm{f(L)} \leq c_{\sigma} \norm{f}_{L^\infty(\Sigma_{\sigma}^0)}.
\]

One can show that $L$ has a bounded holomorphic functional calculus on $H$ if and only if the following quadratic estimates are satisfied:\\
For some (all) $\sigma \in (\omega,\pi)$ and some $\psi \in \Psi(\Sigma_{\sigma}^0)\setminus \{0\}$ there exists some $C>0$ such that for all $x \in H$ 
\begin{align} \label{square-functions}
		C^{-1} \norm{x}^2 \leq \int_0^{\infty} \norm{\psi(tL)x}^2 \,\frac{dt}{t} \leq C \norm{x}^2.
\end{align}

Moreover, if $\psi,\tilde{\psi} \in \Psi(\Sigma_{\sigma}^0)\setminus\{0\}$ are chosen to satisfy $\int_0^{\infty} \psi(t) \tilde{\psi}(t) \, \frac{dt}{t}=1$, then the functional calculus of $L$ on $H$ yields the following \emph{Calder\'{o}n reproducing formula}: For every $f \in H$ 
\begin{align*}
		\int_0^{\infty} \psi(t^{2m}L) \tilde{\psi}(t^{2m}L) f \,\frac{dt}{t} = f	\qquad \text{in}\ H.
\end{align*}
Observe that for given $\psi \in \Psi(\Sigma_{\sigma}^0)\setminus\{0\}$ and given $\alpha,\beta>0$, one can always find a function $\tilde{\psi}  \in \Psi_{\alpha,\beta}(\Sigma_{\sigma}^0)\setminus\{0\}$ such that $\int_0^{\infty} \psi(t) \tilde{\psi}(t) \, \frac{dt}{t}=1$.

\subsection{Tent spaces and Carleson measures}

We recall the most important definitions and properties of tent spaces and Carleson measures.
For proofs of the results, we refer to \cite{CoifmanMeyerStein}. As mentioned in \cite{Stein2}, Chapter II, the proofs, given there in the case of the Euclidean space $\R^n$, carry over to spaces of homogeneous type. \\
For any $x \in X$ , we denote by $\Gamma(x)$ the \emph{cone} of aperture $1$ with vertex $x$, namely
\[
		\Gamma(x):=\{(y,t) \in X \times (0,\infty) \,:\, d(y,x) < t\}.
\]
If $O$ is an open subset of $X$, then the \emph{tent} over $O$, denoted by $\hat{O}$, is defined as
\begin{align*}
		\hat{O}:= \{(x,t) \in X \times (0,\infty) \,:\, \dist(x,O^c) \geq t\}.
\end{align*}

\begin{Def}
For any measurable function $F$ on $X \times (0,\infty)$, the conical square function $\scrA F$ is defined by 
\[
 	\scrA F(x) := \left(\iint_{\Gamma(x)} \abs{F(y,t)}^2 \; \frac{d\mu(y)}{V(x,t)}\frac{dt}{t} \right)^{1/2}, \qquad x \in X,
\]
and the Carleson function $\scrC F$ by
\[
 	\scrC F(x):= \sup_{B \,:\, x \in B} \left(\frac{1}{V(B)} \iint_{\hat{B}} \abs{F(y,t)}^2 \frac{d\mu(y)dt}{t}\right)^{1/2}, \qquad x \in X,
\]
where the supremum is taken over  all balls $B$ in $X$ that contain $x$. \\
For $0<p<\infty$, the tent spaces on $X \times (0,\infty)$ are defined by
\[
 	T^p(X):= \{F: X \times (0,\infty) \to \C  \ \text{measurable} \,;\, \norm{F}_{T^p(X)}:=\norm{\scrA F}_{L^p(X)} < \infty \}.
\]
The tent space $T^{\infty}(X)$ is defined by
\[
 	T^{\infty}(X) := \{F: X \times (0,\infty) \to \C \ \text{measurable} \,;\, \norm{F}_{T^{\infty}(X)} := \norm{\scrC F}_{L^{\infty}(X)} < \infty\}.
\]
\end{Def}

When $p \in [1,\infty]$, the space $(T^p(X), \norm{\,.\,}_{T^p(X)})$ is a Banach space. Moreover, one can show the following duality results.

\begin{Theorem} \label{CMS}
(i) Let $1<p<\infty$ and $\frac{1}{p}+\frac{1}{p'}=1$. There exists a constant $C>0$ such that for all $F \in T^p(X)$ and all $G \in T^{p'}(X)$ 
\[
			\iint_{X \times (0,\infty)} \abs{F(x,t) G(x,t)} \,\frac{d\mu(x)dt}{t}
			\leq C \int_{X} \scrA(F)(x) \scrA(G)(x) \,d\mu(x).
\]
Further, there exists a constant $C>0$ such that for all $F \in T^1(X)$ and all $G \in T^{\infty}(X)$ 
\[
 		\iint_{X \times (0,\infty)} \abs{F(x,t) G(x,t)} \,\frac{d\mu(x)dt}{t}
			\leq C \int_{X} \scrA(F)(x) \scrC(G)(x) \,d\mu(x).
\]
(ii) 
The pairing 
	\[
				\skp{F,G} \mapsto \iint_{X \times (0,\infty)} F(x,t) G(x,t) \,\frac{d\mu(x)dt}{t}
	\]
realizes $T^{p'}(X)$ as equivalent to the dual of $T^p(X)$ if $1<p<\infty$ and $\frac{1}{p}+\frac{1}{p'}=1$, and realizes  $T^{\infty}(X)$ as equivalent to the dual of $T^1(X)$.
\end{Theorem}

We finally state the definition of non-tangential maximal functions and Carleson measures and the connection between both.

\begin{Def}
For any measurable function $F$ on $X \times (0,\infty)$, the \emph{non-tangential maximal function} $F^{\ast}$ is defined by 
\begin{equation} \label{classical-nontang-maxfct}
 	F^{\ast}(x):= \sup_{(y,t) \in \Gamma(x)} \abs{F(y,t)}, \qquad x \in X.
\end{equation}
The space $\calN$ is defined by
$
		\calN:= \{ F: X \times (0,\infty) \to \C  \ \text{measurable} \,;\, \norm{F}_{\calN} := \norm{F^{\ast}}_{L^1(X)} < \infty \}.
$
A \emph{Carleson measure} is a Borel measure $\nu$ on $X \times (0,\infty)$ such that 
\[
 	\norm{\nu}_{\calC} \, := \sup_B  \frac{1}{V(B)} \iint_{\hat{B}} \,\abs{d\nu} < \infty,
\]
where the supremum  is taken over all balls $B$ in $X$. We define $\calC$ to be the space of all Carleson measures.
\end{Def}

The spaces $(\calN,\norm{\,.\,}_{\calN})$ and $(\calC,\norm{\,.\,}_{\calC})$ are Banach spaces. 
 Observe that for $F \in T^{\infty}(X)$ 
\begin{align} \label{tent-carleson}
 	\norm{F}_{T^{\infty}(X)}^2 = \norm{\scrC F}^2_{L^{\infty}(X)} = \norm{\abs{F(y,t)}^2 \,\frac{d\mu(y)dt}{t}}_{\calC}.
\end{align}

\begin{Theorem} \label{CarlesonDuality}
 If $F \in \calN$ and $\nu \in \calC$, then 
\[	
 	\iint_{X \times (0,\infty)} \abs{F(x,t)} \, d\nu(x,t) \leq C \norm{F}_{\calN} \cdot \norm{\nu}_{\calC}.
\]
\end{Theorem}

For applications, we also need the following corollary. 

\begin{Prop} \label{Carleson-Tent-Cor}
Let $2<p<\infty$. Let $F$ be a measurable function on $X \times (0,\infty)$ with $F^{\ast} \in L^p(X)$ and let $G \in T^{\infty}(X)$. Then 
\begin{align*}
		\norm{\scrC(F \cdot G)}_{L^p(X)} \leq C \norm{F^{\ast}}_{L^p(X)} \norm{\scrC G}_{L^{\infty}(X)},
\end{align*}
with a constant $C>0$ independent of $F$ and $G$.
\end{Prop}

\subsection{Assumptions on the operator}
\label{sect-assumpOperator}

We fix our assumptions on the operator $L$. 
Unless otherwise specified, we will assume the following.

\begin{enumerate}[label= \textbf{(H\arabic*)},ref=H\theenumi]
\item \label{H1} 
The operator $L$ is an injective, sectorial operator in $L^2(X)$ of angle $\omega$, where $0 \leq \omega < \pi/2$. Further, $L$ has a bounded $H^{\infty}(\Sigma_\sigma^0)$-functional calculus for some (all) $\omega < \sigma < \pi$.
\item \label{H2}
	The operator $L$ generates an analytic semigroup $\{e^{-tL}\}_{t>0}$ satisfying Davies-Gaffney condition. 
	That is, there exist constants $C,c>0$ such that for arbitrary open subsets $E,F \subseteq X$ 
\begin{equation} \label{DaviesGaffney}
		\norm{e^{-tL} f}_{L^2(F)}
			\leq C \exp\left[ - \left(\frac{\dist(E,F)^{2m}}{ct}\right)^{\frac{1}{2m-1}}\right] \norm{f}_{L^2(E)} 
\end{equation}
for every $t>0$ and every $f \in L^2(X)$ with $\supp f \subseteq E$.
\end{enumerate}

For the theory of Hardy and BMO spaces associated to $L$, these two assumptions will be enough.
In order to show $L^2(X)$-boundedness of certain paraproducts, we need one additional assumption. Henceforth, we will explicitly mention whenever we take into account the following assumption.
\begin{enumerate}[label= \textbf{(H\arabic*)},ref=H\theenumi]
\setcounter{enumi}{2}
\item \label{H3} 
	The semigroup $\{e^{-tL}\}_{t>0}$ satisfies an $L^{\tilde{p}}-L^2$ off-diagonal estimate for some $\tilde{p} \in (1,2)$ and an $L^2-L^{\tilde{q}}$ off-diagonal estimate for some $\tilde{q} \in (2,\infty)$, i.e. there exists a constant $C>0$ and some $\eps>0$ such that for every $t>0$, every $j \in \N_0$ and for an arbitrary ball $B$ in $X$ with radius $r=t^{1/2m}$ 
\begin{align} \label{Lp-L2-estimate}
 	\norm{e^{-tL} \Eins_{S_j(B)} f}_{L^2(B)} \leq C 2^{-j(\frac{n}{\tilde{p}}+\eps)} V(B)^{\frac{1}{2}-\frac{1}{\tilde{p}}} \norm{f}_{L^{\tilde{p}}(S_j(B))}
\end{align}
and
\begin{align} \label{L2-Lq-estimate}
 	\norm{e^{-tL} \Eins_B g}_{L^{\tilde{q}}(S_j(B))} \leq C 2^{-j(\frac{n}{\tilde{q}'}+\eps)} V(B)^{\frac{1}{\tilde{q}}-\frac{1}{2}} \norm{g}_{L^2(B)}	
\end{align}
for all $f \in L^{\tilde{p}}(X)$ and all $g \in L^2(X)$.
Here, $\tilde{q}'$ is the conjugate exponent of $\tilde{q}$ defined by $\frac{1}{\tilde{q}}+\frac{1}{\tilde{q}'}=1$.
\end{enumerate}
Observe that \eqref{L2-Lq-estimate} is just the dual estimate of \eqref{Lp-L2-estimate}. That is, if $L$ satisfies \eqref{L2-Lq-estimate} with exponent $\tilde{q}$, then $L^{\ast}$ satisfies \eqref{Lp-L2-estimate} with exponent $\tilde{q}'$ and vice versa.\\

One can show that the Davies-Gaffney estimates imply $L^2$ off-diagonal estimates for more general operator families associated to $L$. The proof of \cite{HofmannMayborodaMcIntosh}, Lemma 2.28, carries over with only minor changes to our more general setting.

\begin{Prop}  \label{H-inf-offdiag}
 Let $L$ satisfy \eqref{H1} and \eqref{H2}. Let $\sigma \in (\omega, \frac{\pi}{2}), \ \psi \in \Psi_{\alpha,\beta}(\Sigma_{\sigma}^0)$ for some $\alpha, \beta >0$ and $\varphi \in H^{\infty}(\Sigma_{\sigma}^0)$. Then the family of operators $\{\psi(tL)\varphi(L)\}_{t>0}$ satisfies $L^2$ off-diagonal estimates of order $\alpha$, with the constant controlled by $\norm{\varphi}_{L^{\infty}(\Sigma_{\sigma}^0)}$.
That is, there exists a constant $C>0$ such that for arbitrary open sets $E,F \subseteq X$
\begin{align*}
 	\norm{\psi(tL)\varphi(L)f}_{L^2(F)} \leq C \norm{\varphi}_{L^{\infty}(\Sigma_{\sigma}^0)} \left(1+\frac{\dist(E,F)^{2m}}{t}\right)^{-\alpha} \norm{f}_{L^2(E)}
\end{align*}
for every $t>0$ and every $f \in L^2(X)$ supported in $E$.
\end{Prop}

We end the section with an observation on conservation properties of the semigroup. 

\begin{Lemma} \label{psi-remark}
Let $L$ satisfy \eqref{H1}, \eqref{H2} and let $\sigma \in (\omega, \frac{\pi}{2})$.\\
(i) Let $\gamma>\frac{n}{4m}$. For every ball $B \subseteq X$ there exists some constant $C_B>0$ such that for all $t>0$ 
\begin{align*}
		\norm{e^{-tL^{\ast}}}_{L^2(B)\to L^1(X\setminus 4B)} \leq C_B t^{\gamma}.
\end{align*}
In particular, one can define via duality $e^{-tL}$ as an operator from $L^{\infty}(X)$ to $L^2_{\loc}(X)$.\\
(ii) Let $\alpha>0$, $\beta>\frac{n}{4m}$ and $\psi \in \Psi_{\beta,\alpha}(\Sigma_{\sigma}^0)$. Moreover, let $b \in L^{\infty}(X)$. 
If for every $t>0$ 
\[
 	e^{-tL}(b)=b  \qquad \text{in} \ L^2_{\loc}(X),
\]
then for every $t>0$
\[
 	\psi(tL)(b)=0 \qquad \text{in} \ L^2_{\loc}(X).
\]
\end{Lemma}

\begin{Proof}
(i) Let $f \in L^2(X)$ with $\supp f \subseteq B$. Due to the Cauchy-Schwarz inequality, \eqref{H2} and the doubling condition \eqref{doublingProperty2} there holds
\begin{align*}
		\norm{e^{-tL^{\ast}}f}_{L^1(X\setminus 4B)} 
				&\leq \sum_{j=1}^{\infty} V(2^jB)^{1/2} \norm{e^{-tL^{\ast}}f}_{L^2(S_j(B))} \\
				& \lesssim \sum_{j=1}^{\infty} V(2^jB)^{1/2} \exp\left(-\frac{\dist(B,S_j(B))^{2m}}{t}\right) \norm{f}_{L^2(B)} \\
				& \lesssim V(B)^{1/2} \sum_{j=1}^{\infty} 2^{jn/2} \left(\frac{t}{(2^jr_B)^{2m}}\right)^{\gamma}  \norm{f}_{L^2(B)} 
				 \leq C_B t^{\gamma} \norm{f}_{L^2(B)},
\end{align*}
where in the last step we used the assumption $\gamma>\frac{n}{4m}$.\\
(ii) Let $\gamma \in (\frac{n}{4m},\beta)$. Moreover, let $\omega<\theta<\sigma<\frac{\pi}{2}$ and $\lambda \in \partial \Sigma_{\theta}^0$. 
According to (i), the integral
\begin{align*}
		\int_0^{\infty} e^{-\lambda t} e^{-tL^{\ast}} \,dt
\end{align*}
converges strongly as an operator from $L^2(B)$ to $L^1(X\setminus 4B)$ with the operator norm bounded by a constant times $\abs{\lambda}^{-\gamma-1}$.
This also implies that $\nnorm{\psi(\lambda)(\lambda+L^{\ast})^{-1}}_{L^2(B) \to L^1(X\setminus 4B)} \lesssim \abbs{\psi(\lambda)}\abs{\lambda}^{-\gamma-1}$ and the integral
\begin{align*}
		\frac{1}{2\pi i} \int_{\partial \Sigma_{\theta}^0} \psi(\lambda) (\lambda+L^{\ast})^{-1} \,d\lambda,
\end{align*}
converges strongly as an operator from $L^2(B)$ to $L^1(X)$, since $\beta>\gamma$.
The assumption $e^{-tL}(b)=b$ then yields for every $f \in L^2(B)$
\begin{align*}
	\skp{b,(\lambda+L^{\ast})^{-1}f} = \skp{b,\int_0^{\infty} e^{-\lambda t} e^{-tL^{\ast}} f \,dt}
				=\int_0^{\infty} e^{-\lambda t} \skp{e^{-tL}(b),f} \,dt
				=\frac{1}{\lambda} \skp{b,f}.
\end{align*}
We finally obtain for $\psi(L)(b)$ the equality
\begin{align*}
		\skp{\psi(L)(b),f} 
				 = \skp{b,\psi(L^{\ast})f}
				 = \frac{1}{2\pi i} \int_{\partial \Sigma_{\theta}^0} \psi(\lambda) \skp{b,(\lambda+L^{\ast})^{-1}f} \,d\lambda
				 =  \frac{1}{2\pi i} \int_{\partial \Sigma_{\theta}^0} \frac{\psi(\lambda)}{\lambda} \,d\lambda \; \skp{b,f} 
				 = 0,
\end{align*}
where the last step is due to an extension of Cauchy's theorem and the assumption $\psi \in \Psi(\Sigma_{\sigma}^0)$.
\end{Proof}

\section{Hardy and BMO spaces associated to operators revisited}

In the following, we will always assume that the operator $L$ satisfies the assumptions \eqref{H1} and \eqref{H2} and that $\sigma \in (\omega,\frac{\pi}{2})$.
We denote by $\calD(S)$ the domain, by $\calR(S)$ the range of an unbounded operator $S$, and by $S^k$ the $k$-fold composition of $S$ with itself, in the sense of unbounded operators.\\
We summarize the most important facts about Hardy and BMO spaces associated to $L$. For more details and proofs of the results, we refer to  \cite{HofmannMayboroda}, \cite{HofmannMayborodaMcIntosh}, \cite{HLMMY} and \cite{DuongLi}. 
The proofs given there carry over with only minor changes to our more general setting. 
In addition, we generalize a Calder\'on reproducing formula for elements of $H^1_L(X)$ and $BMO_{L^{\ast}}(X)$ and a Carleson measure estimate. Both results have their origin in \cite{HofmannMayboroda}.

\subsection{The spaces $H^p_L(X)$ and $BMO_L(X)$}

Let $\psi \in \Psi(\Sigma_{\sigma}^0)\setminus \{0\}$ and consider for every $f \in L^2(X)$ the square function $\scrA Q_{\psi,L}f$ associated to $L$, namely
 \[
 			\scrA Q_{\psi,L}(f)(x) = \left(\iint_{\Gamma(x)} \abs{\psi(t^{2m}L) f(y)}^2 \, \frac{d\mu(y)}{V(x,t)}\frac{dt}{t} \right)^{1/2}, \quad x \in X.
 	\]

\begin{Def} \label{Def-Hp}
(i) Let $1 \leq p \leq 2$ and let $\psi_0 \in \Psi(\Sigma_{\sigma}^0)$ be defined by $\psi_0(z):=ze^{-z}$. Define $H^p_L(X)$  to be the completion of the space
\begin{equation} \label{Def-bH1}
		\bH^p_L(X):= \{f \in L^2(X) \,:\, \scrA Q_{\psi_0,L}f \in L^p(X)\},
\end{equation}
 with respect to the norm
$
 			\norm{f}_{H^p_{\psi_0,L}(X)} := \norm{\scrA Q_{\psi_0,L}f}_{L^p(X)}. 
$\\
(ii) Let $2<p<\infty$. Define
$
			H^p_L(X) :=(H^{p'}_{L^{\ast}}(X))',
$
where $\frac{1}{p}+\frac{1}{p'}=1$ and $L^{\ast}$ is the adjoint operator of $L$.
\end{Def}

Observe that $\norm{\scrA Q_{\psi,L}f}_{L^p(X)} = \norm{Q_{\psi,L}f}_{T^p(X)}$ for $Q_{\psi,L}f(x,t):=\psi(t^{2m}L)f(x)$. Moreover, there holds $H^2_L(X)=L^2(X)$. \\

In both cases, for $p \leq 2$ and for $p>2$, there is a characterization of $H^p_L(X)$ by general square functions constructed via functions $\psi \in \Psi(\Sigma_{\sigma}^0) \setminus \{0\}$ with a certain decay at infinity and at zero, respectively. 
For a proof, we refer to Corollary 4.21 of \cite{HofmannMayborodaMcIntosh}.

\begin{Theorem} \label{Charact-Hp}
	Let $\alpha>0$ and $\beta>\frac{n}{4m}$. Further, let either $1\leq p \leq 2$ and $\psi \in \Psi_{\alpha,\beta}(\Sigma_{\sigma}^0) \setminus \{0\}$ or $2 \leq p < \infty$ and $\psi \in \Psi_{\beta,\alpha}(\Sigma_{\sigma}^0) \setminus \{0\}$.
	Define $H^p_{\psi,L}(X)$ to be the completion of the space
	\[
				\bH^p_{\psi,L}(X) :=\{f \in L^2(X) \,:\,  \scrA Q_{\psi,L}f \in L^p(X)\},
	\]
	 with respect to the norm
$
 			\norm{f}_{H^p_{\psi,L}(X)} := \norm{\scrA Q_{\psi,L}f}_{L^p(X)}. 
$
Then $H^p_L(X)=H^p_{\psi,L}(X)$, with equivalence of norms.
\end{Theorem}

There also exists a molecular characterization of $H^1_L(X)$. 
We begin with a definition of molecules associated to $L$.

 \begin{Def} \label{DefMolecules}
  Let $M \in \N$ and $\eps>0$. A function $m \in L^2(X)$ is called a $(1,2,M,\eps)$-molecule associated to $L$ if there exists a function $b \in \calD(L^M)$ and a ball $B$ in $X$ with radius $r_B>0$ such that \\
(i) $m=L^M b$;\\
(ii) For every $k=0,1,2,\ldots,M$ and all $j \in \N_0$ 
 \[
  	\norm{(r_B^{2m}L)^k b}_{L^2(S_j(B))} \leq r_B^{2mM} 2^{-j\eps} V(2^jB)^{-1/2}.
 \]
 \end{Def}

The molecular Hardy spaces associated to $L$ are then defined as follows.

 \begin{Def}
 Given $M \in \N$, $\eps>0$ and $f \in L^1(X)$, one says that $f=\sum_j \lambda_j m_j$ is a \emph{molecular $(1,2,M,\eps)$-representation of $f$} if $\sum_{j=0}^{\infty} \abs{\lambda_j} < \infty$, each $m_j$ is a $(1,2,M,\eps)$-molecule, and the sum converges in $L^2(X)$.\\
Let $\eps>0$ be fixed. Set
 \[
  	\bH_{L,mol,M}^1(X) := \{f \in L^1(X)\,:\, f \;\text{has a}\; (1,2,M,\eps)\text{-representation}\}
 \]
 with the norm given by
 \begin{align*}
   \norm{f}_{H^1_{L,mol,M}(X)} := 
 	\inf \left\{\sum_{j=0}^{\infty} \abs{\lambda_j} \, : \, f= \sum_{j=0}^{\infty} \lambda_j m_j \, \; \text{is a} \; (1,2,M,\eps)\text{-representation} \right\}.
 \end{align*}

  The space $H^1_{L,mol,M}(X)$ is defined to be the completion of $\bH^1_{L,mol,M}(X)$ with respect to the norm $\norm{\,.\,}_{H^1_{L,mol,M}(X)}$ above.
 \end{Def}

One can show the following equivalence. For a proof, we refer to \cite{DuongLi}, Theorem 3.12.

 \begin{Theorem} \label{CharacterizationHardy}
   Suppose that $M \in \N$ with $M>\frac{n}{4m}$  Then
   $
   				H^1_{L,mol,M}(X) = H^1_L(X)
   	$
   	with equivalence of norms.
 \end{Theorem}

Next, let us define the space $BMO_L(X)$. Let us fix some element $x_0 \in X$ that will henceforth be called $0$. The ball $B_0:=B(0,1)$ will then be referred to as \emph{unit ball}.
One first defines a space $\calE_M(L)$ in such a way that for every $f \in \calE_M(L)$ there holds $(I-e^{r_B^{2m}L})^M f \in L^2_{\loc}(X)$, and therefore the expression in \eqref{defBMO-1} is well-defined.

\begin{Def}
Let $\eps>0$, $M \in \N$ and let $\phi \in \calR(L^M) \subseteq L^2(X)$ with $\phi = L^M \nu$ for some $\nu \in \calD(L^M)$. Introduce the norm
\[
 	\norm{\phi}_{\molMeps} := \sup_{j \geq 0} \left[2^{j\eps} V(2^jB_0)^{1/2} \sum_{k=0}^M \norm{L^k \nu}_{L^2(S_j(B_0))}\right],
\]
where $B_0$ is the unit ball centered at $0$ with radius 1, and set
\[
 	\molMeps:=\{\phi \in \calR(L^M) \,:\, \norm{\phi}_{\molMeps} < \infty\}.
\]
One denotes by $(\molMeps)'$ the dual of $\molMeps$. 
For any $M \in \N$, let $\calE_M(L)$ be defined by
\[
 	\calE_M(L):= \bigcap_{\eps>0} (\molMdeps)'.
\]
\end{Def}

\begin{Remark}
Let $M \in \N$ and $\eps>0$. 
Then for every $f \in (\molMdeps)'$ and every $t>0$, one can via duality define $(I-e^{-t^{2m}L})^M f$ and $(I-(I+t^{2m}L)^{-1})^Mf$ as elements of $L^2_{\loc}(X)$.
\end{Remark}

\begin{Def}
 Let  $M \in \N$. An element $f \in \calE_M(L)$ is said to belong to $BMO_{L,M}(X)$ if
\begin{equation} \label{defBMO-1}
 	\norm{f}_{BMO_{L,M}(X)} := 
		\sup_{B \subseteq X} \left(\frac{1}{V(B)} \int_B \abs{(I-e^{-r_B^{2m}L})^M f(x)}^2 d\mu(x) \right)^{1/2} < \infty,
\end{equation}
where the supremum is taken over all balls $B$ in $X$.
\end{Def}

One can then show the following duality result. For a proof, we refer to \cite{DuongLi}, Theorem 3.23 and 3.24. 

 \begin{Theorem} \label{Duality}
Let $M>\frac{n}{4m}$. Then 
$
 	(H^1_L(X))' = BMO_{L^{\ast},M}(X).
$
 \end{Theorem}

 In particular, the theorem yields that the definition of $BMO_{L,M}(X)$ is independent of the choice of $M>\frac{n}{4m}$. This leads to the following definition.
\begin{Def} \label{DefBMO-uniform}
The space $BMO_L(X)$ is defined by
$
		BMO_L(X):=BMO_{L,M}(X),
$
where $M \in \N$ with $M>\frac{n}{4m}$.
\end{Def}

\subsection{Interpolation of Hardy spaces}

The spaces $H^p_L(X)$ form a complex interpolation scale.
For a proof, we refer to \cite{HofmannMayborodaMcIntosh}, Lemma 4.24, where the authors reduce the problem to complex interpolation of tent spaces.

\begin{Prop} \label{Interpolation}
 Let $L$ be an operator satisfying \eqref{H1} and \eqref{H2}. Let $1 \leq p_0<p_1<\infty$ and $0<\theta<1$.
 Then
 \begin{align*}
 		[H^{p_0}_L(X),H^{p_1}_L(X)]_{\theta} & = H^p_L(X) \qquad \text{where} \;\, 1/p=(1-\theta)/p_0 + \theta/p_1, \\
 				[H^{p_0}_L(X),BMO_L(X)]_{\theta} &= H^p_L(X) \qquad \text{where} \;\, 1/p=(1-\theta)/p_0.
 \end{align*}
\end{Prop}

The next result is a slight generalization of \cite{HofmannMayboroda}, Theorem 3.2 and complements \cite{BlunckKunstmann}, Theorem 1.1.  

\begin{Prop} \label{H1-bddness}
Let $M \in \N$, $M>\frac{n}{4m}$. Assume that $T$ is a linear or a non-negative sublinear operator defined on $L^2(X)$ such that
$
 	T: L^2(X) \to L^2(X)
$
is bounded and $T$ satisfies the following weak off-diagonal estimates:\\
There exists some $\gamma > \frac{n}{2m}$ and a constant $C>0$ such that for every $t>0$, arbitrary balls $B_1,B_2 \in X$ with radius $r=t^{1/2m}$ and every $f \in L^2(X)$ supported in $B_1$
\begin{align} \label{H1-eq1}
 	\norm{T(I-e^{-tL})^M(f)}_{L^2(B_2)} &\leq C_T \left(1+\frac{\dist(B_1,B_2)^{2m}}{t}\right)^{-\gamma} \norm{f}_{L^2(B_1)},  \\ \label{H1-eq2}
	\norm{T(tLe^{-tL})^M(f)}_{L^2(B_2)} &\leq C_T \left(1+\frac{\dist(B_1,B_2)^{2m}}{t}\right)^{-\gamma} \norm{f}_{L^2(B_1)}.
\end{align}
Then
$
 	T: H^1_L(X) \to L^1(X)
$
is bounded and there exists some $C>0$, independent of $C_T$, such that for all  $f \in H^1_L(X)$ 
\begin{align*}
		\norm{Tf}_{L^1(X)} \leq C C_T \norm{f}_{H^1_L(X)}.
\end{align*}
\end{Prop}

\begin{Remark}
 If \eqref{H1-eq1} and \eqref{H1-eq2} are satisfied for arbitrary open sets $E,F \subseteq X$, one only requires a decay of order $\gamma>\frac{n}{4m}$.
\end{Remark}

A sufficient condition and a detailed proof for the equivalence of $H^p_L(X)$ and $L^p(X)$ is given in \cite{Uhl}, Theorem 4.19.
We refer the reader to a comparison with assumption \eqref{H3}.

\begin{Prop} \label{Hp-equiv}
Let $L$ satisfy \eqref{H1} and \eqref{H2}.	If for some $p_0 \in [1,2)$, there exist constants $C,c>0$ such that for all $x,y \in X$ and all $t>0$
\begin{align*}
		& \norm{\Eins_{B(x,t^{1/2m})} e^{-tL} \Eins_{B(y,t^{1/2m})}}_{L^{p_0}(X) \to L^{p'_0}(X)} \\
				& \qquad \qquad \qquad \qquad \leq C V(x,t^{1/2m})^{-(\frac{1}{p_0} - \frac{1}{p'_0})} \exp\left(-\left(\frac{d(x,y)^{2m}}{ct}\right)^{\frac{1}{2m-1}}\right),
\end{align*}
then
\begin{align*}
		H^p_L(X) = L^p(X), \qquad p_0<p < p'_0.
\end{align*}	
\end{Prop}

For further relationships between $H^p_L(X)$ and $L^p(X)$ in the case of second order elliptic operators in divergence form, we refer to \cite{HofmannMayborodaMcIntosh}, Proposition 9.1.

\subsection{A Calder\'{o}n reproducing formula and Carleson measures} 
As shown in \cite{HofmannMayboroda}, Lemma 8.4, it is possible to generalize the Calder\'{o}n reproducing formula, originally given on $L^2(X)$ via functional calculus, to functions $f \in BMO_{L^{\ast},M}(X)$ and functions $g \in H^1_L(X)$, that can be represented as a finite linear combination of molecules. 
Compared with \cite{HofmannMayboroda}, Lemma 8.4, we state a more general version of the lemma, that allows for a greater freedom in the choice of functions $\psi, \, \tilde{\psi} \in \Psi(\Sigma_{\sigma}^0)$.

 \begin{Lemma} \label{ReproducingFormula}
  Let $M \in \N$ and suppose that $f \in \calE_M(L^{\ast})$ satisfies the ``controlled growth estimate''
 \begin{equation} \label{growthBound}
  	\int_{X} \frac{\abs{(I-(I+L^{\ast})^{-1})^M f(x)}^2}{(1+d(x,0))^{\eps_1}V(0,1+d(x,0))} \, d\mu(x) < \infty
 \end{equation}
 for some $\eps_1>0$.
Let $\psi \in \Psi_{\beta_1,\alpha_1}(\Sigma_\sigma^0)\setminus\{0\}$ and $\tilde{\psi} \in \Psi_{\beta_2,\alpha_2}(\Sigma_\sigma^0)\setminus\{0\}$ for some constants $\alpha_1,\alpha_2,\beta_1,\beta_2 >0$, with $\beta_1 + \beta_2 > \frac{n+\eps_1}{4m}$ and $\int_0^{\infty} \psi(t) \tilde{\psi}(t) \frac{dt}{t}=1$.
 Then for every $g \in H^1_L(X)$ that can be represented as a finite linear combination of $(1,2,M',\eps)$-molecules, with $\eps> \frac{\eps_1}{2}$, $M'-M> \frac{n+\eps_1}{4m}$ and $\alpha_1 + \alpha_2 >M'$, we have
 \[
  	\skp{f,g} 
 		=  \lim_{\substack{\delta \to 0\\ R \to \infty}} \int_{\delta}^R \int_X \psi(t^{2m}L^{\ast}) f(x) \overline{\tilde{\psi}(t^{2m}L) g(x)} \;\frac{d\mu(x)dt}{t}.
 \]
 \end{Lemma}

 \begin{Remark} \label{Remark-BMO}
If $f \in BMO_{L^{\ast},M}(X)$, then condition \eqref{growthBound} is fulfilled for every $\eps_1>0$.
 \end{Remark}

The proof works in most parts analogously to the one of \cite{HofmannMayboroda}. We need one lemma in addition, which gives us a primitive of a function $\psi \in \Psi(\Sigma_\sigma^0)$.

 \begin{Lemma} \label{primitive}
  Let $\sigma \in (0,\pi)$, $\alpha,\beta>0$ and $\psi \in \Psi_{\beta,\alpha}(\Sigma_{\sigma}^0)\setminus\{0\}$. Then for every $l \in \N$ with $l \geq \alpha$  there exists a function $\varphi \in \Psi_{\beta,\alpha}(\Sigma_{\sigma}^0)$ and some $\gamma \in \C$ such that
 \[
  	\psi(z) = z \varphi'(z) + \gamma \frac{z}{(1+z)^{l+1}}, \qquad z \in \Sigma_{\sigma}^0.
 \]
 \end{Lemma}
 
 \begin{Proof}
  Let us define a function $G$ on $\Sigma_{\sigma}^0$ by setting
\[
 	G(z):= \int_{\gamma_z} \frac{\psi(\zeta)}{\zeta} \,d\zeta, \qquad z \in \Sigma_{\sigma}^0,
\]
where $\gamma_z(t):= te^{i \arg z}, \; t \geq \abs{z}$, is the parametrization of the half-ray with angle $\arg z$ starting at $z$.
By assumption there holds $\frac{\psi(\zeta)}{\zeta} = \calO(\abs{\zeta}^{-\alpha-1})$ for $\abs{\zeta} \to \infty$ and consequently,
$G(z)=\calO(\abs{z}^{-\alpha})$ for $\abs{z} \to \infty$.
By definition of $G$, we further have 
\[
 	zG'(z) = \psi(z), \qquad z \in \Sigma_{\sigma}^0.
\]
To get the desired behaviour at $0$, one has to do a little more work. We know by assumption that $\frac{\psi(z)}{z} = \calO(\abs{z}^{\beta-1})$ for $\abs{z} \to 0$ and, since $\beta>0$, the integral
\begin{align} \label{primitive-int}
 	\int_{\Gamma_{\theta}} \frac{\psi(\zeta)}{\zeta} \,d\zeta
\end{align}
converges for every $\theta \in (-\sigma,\sigma)$, where $\Gamma_{\theta}(t):=te^{i\theta}, \; 0<t<\infty$. Using the same arguments as in \cite{LevicoNotes}, Remark 9.3, one can show that due to Cauchy's theorem, the integral in \eqref{primitive-int} is independent of the angle $\theta \in (-\sigma,\sigma)$. 
Therefore, let us set $c:=\int_{\Gamma_{\theta}} \frac{\psi(\zeta)}{\zeta} \,d\zeta$ for any $\theta \in (-\sigma,\sigma)$.
We then obtain
\[
 	c-G(z) = \int_{\widetilde{\gamma}_z} \frac{\psi(\zeta)}{\zeta} \,d\zeta, \qquad z \in \Sigma_{\sigma}^0,
\]
where $\widetilde{\gamma}_z (t) := te^{i\arg z}$, $0 < t \leq \abs{z}$, is the parametrization of the half-ray with angle $\arg z$ starting at $0$ and ending at $z$. 
From the assumption $\frac{\psi(\zeta)}{\zeta} = \calO(\abs{z}^{\beta-1})$ for $\abs{z} \to 0$ we now get that 
$c - G(z) = \calO(\abs{z}^{\beta})$ for $\abs{z} \to 0$. 
Therefore, by defining for a given $l \in \N$ with $l \geq \alpha$
\begin{align*}
 	\varphi(z):= G(z) - \frac{c}{(1+z)^l}, \qquad z \in \Sigma_{\sigma}^0,
\end{align*}
we obtain the following:
By construction there holds $\varphi(z) = \calO(\abs{z}^{\beta})$ for $\abs{z} \to 0$ and $\varphi(z) = \calO(\abs{z}^{-\alpha})$ for $\abs{z} \to \infty$.
In addition, a simple calculation shows that
\[
 	\psi(z) = z G'(z) = z\varphi'(z) - \frac{lcz}{(1+z)^{l+1}},
\]
which concludes the proof with $\gamma=-lc$.
 \end{Proof}

%
%

The relation of elements of $BMO_L(X)$ and Carleson measures can be described as follows. 

  \begin{Prop} \label{LemmaCarlesonMeasure}
Let $M \in \N, \; M > \frac{n}{4m}$. Further, let $\alpha>0$, $\beta> \frac{n}{4m}$ and $\psi \in \Psi_{\beta,\alpha}(\Sigma_{\sigma}^0)\setminus\{0\}$. Then the operator
 \[
  	f \mapsto \psi(t^{2m}L)f
 \]
 maps $BMO_{L}(X) \to T^{\infty}(X)$, i.e. for every $f \in BMO_{L}(X)$ is
 \begin{equation} \label{CarlesonMeasure}
  	\nu_{\psi,f} := \abs{\psi(t^{2m}L)f(y)}^2 \frac{d\mu(y)\, dt}{t}
 \end{equation}
 a Carleson measure and there exists a constant $C_{\psi}>0$ such that for all $f \in BMO_{L}(X)$
 \[
  	\norm{\nu_{\psi,f}}_{\calC} \leq C_{\psi} \norm{f}^2_{BMO_{L}(X)}.
 \]
Conversely, if $f \in \calE_M(L)$ satisfies the controlled growth bound \eqref{growthBound} (with $L$ in place of $L^{\ast}$) for some $\eps_1>0$, and if $\nu_{\psi,f}$ defined in \eqref{CarlesonMeasure} is a Carleson measure, then $f \in BMO_L(X)$ and
 \[
  	\norm{f}_{BMO_L(X)}^2 \leq C \norm{\nu_{\psi,f}}_{\calC}.
 \]
 \end{Prop}

For a special choice of $\psi$, namely $\psi(z)=z^Me^{-z}$, the result is due to \cite{HofmannMayboroda}, Theorem 9.1. In the generality as stated above, the first part of the result is due to \cite{HofmannMayborodaMcIntosh}, Proposition 4.13. The second part is new and can be shown by combining the proof of \cite{HofmannMayboroda}, Theorem 9.1 with Lemma \ref{ReproducingFormula}.

 \section{Paraproducts via $H^{\infty}$-functional calculus}
\label{sect-paraproduct}

In this section, we introduce paraproduct operators associated to a sectorial operator $L$ and investigate various properties of those. 

\subsection{Boundedness of paraproducts - via Carleson measures}
\label{paraproduct-L2}

We begin with the study of the following paraproduct operator.

 \begin{Def} \label{DefParaproduct}
 Let $L$ satisfy \eqref{H1}. Assume that $\psi, \tilde{\psi} \in \Psi(\Sigma_{\sigma}^0)\setminus \{0\}$. For $b \in BMO_L(X)$ and $f \in L^2(X)$ we define the paraproduct
 \begin{equation} \label{Paraproduct}
 	\Pi_b(f):= \int_0^{\infty} \tilde{\psi}(t^{2m}L) [\psi(t^{2m}L)b \cdot A_t(e^{-t^{2m}L}f)] \, \frac{dt}{t},
 \end{equation}
 where $A_t$ is the averaging operator defined in \eqref{avOp}.
 \end{Def}
 
For convenience, we do not index $\Pi_b$ with the defining functions $\psi$ and $\tilde{\psi}$. In the context, it will always become clear what the defining functions are.

 \begin{Theorem} \label{paraproduct}
  Assume that $L$ satisfies \eqref{H1} and \eqref{H2}. Let $\alpha>0$, $\beta>\frac{n}{4m}$ and let $\psi \in \Psi_{\beta,\alpha}(\Sigma_{\sigma}^0)\setminus \{0\}$. \\
(i) Let $L$ satisfy in addition \eqref{Lp-L2-estimate} of \eqref{H3} and assume that $\tilde{\psi} \in \Psi(\Sigma_{\sigma}^0)\setminus \{0\}$.
Then the operator $\Pi_b$, defined in \eqref{Paraproduct}, is bounded on $L^2(X)$ for every $b \in BMO_L(X)$, i.e. there exists some constant $C>0$ such that for every $f \in L^2(X)$ and every $b \in BMO_L(X)$
 \[
  	\norm{\Pi_b(f)}_{L^2(X)} \leq C \norm{b}_{BMO_L(X)} \norm{f}_{L^2(X)}.
 \]
(ii) Let $p \in (2,\infty]$ and assume that $\tilde{\psi} \in \Psi_{\alpha,\beta}(\Sigma_{\sigma}^0)\setminus \{0\}$.
 Then the operator $\Pi_b$, initially defined on $L^2(X)$ in \eqref{Paraproduct}, extends for every  $b \in BMO_{L}(X)$ to a bounded operator $\Pi_b: L^p(X) \to H^p_L(X)$. That is, there exists some constant $C>0$ such that for every  $b \in BMO_L(X)$  and every $f \in L^p(X)$ 
\begin{align*}
		\norm{\Pi_b(f)}_{H^p_L(X)} \leq C  \norm{b}_{BMO_{L}(X)} \norm{f}_{L^p(X)}.
\end{align*}
Here, we designate $H^{\infty}_L(X):=BMO_L(X)$.
 \end{Theorem}

The combination of Theorem \ref{paraproduct} and Proposition \ref{Hp-equiv} yields appropriate boundedness results on $L^p(X)$ instead of $H^p_L(X)$.\\

We start the preparations for the proof with the following definition of a modified non-tangential maximal function. 
The modification is required in absence of pointwise estimates. It has its origin in \cite{KenigPipher} and was e.g. recently applied in \cite{HofmannMayboroda}.

 \begin{Def} \label{def-nontang-maxfct}
 	Given an operator $L$ satisfying \eqref{H1} and a function $f \in L^2(X)$ we define the non-tangential maximal operator $\calN_{h,L}$ associated to $L$ via 
 \[
 		\calN_{h,L} f(x):= \sup_{(y,t) \in \Gamma(x)} \left(\frac{1}{V(y,t)} \int_{B(y,t)} \abs{e^{-t^{2m}L} f(z)}^2 d\mu(z) \right)^{1/2}, \quad x \in X.
 \]
 \end{Def}

We can then show the following.

\begin{Lemma} \label{nontang-maxfct}
(i)  Assume that $L$ satisfies \eqref{H1} and \eqref{Lp-L2-estimate} of \eqref{H3}. Then the operator $\calN_{h,L}$ is bounded on $L^2(X)$, i.e. there exists a constant $C>0$ such that for every $f \in L^2(X)$ 
 \[
  	\norm{\calN_{h,L}f}_{L^2(X)} \leq C \norm{f}_{L^2(X)}.
 \]
 (ii)  Assume that $L$ satisfies \eqref{H1} and \eqref{H2}. Then the operator $\calN_{h,L}$ is bounded on $L^p(X)$ for every $p \in (2,\infty]$.
 \end{Lemma}

 \begin{Proof}
(i) We will show a pointwise estimate of $\calN_{h,L} f $ against the uncentered maximal function $\calM_{\tilde{p}} f$, where the index $\tilde{p} \in (1,2)$ comes from assumption \eqref{H3}. \\
 Let  $f\in L^2(X)$ and $x \in X$.  To apply the $L^{\tilde{p}}-L^2$ off-diagonal estimates for the semigroup, we use an annular decomposition of $f$.
  This yields
  \begin{align*}
   	\calN_{h,L} f(x) &= \sup_{(y,t) \in \Gamma(x)} \left( \frac{1}{V(y,t)} \int_{B(y,t)} \abs{e^{-t^{2m}L} f(z)}^2 \,d\mu(z) \right)^{1/2} \\
 				& \leq \sup_{(y,t) \in \Gamma(x)} \sum_{j=0}^{\infty} V(y,t)^{-1/2} \norm{ e^{-t^{2m}L} \Eins_{S_j(B(y,t))} f}_{L^2(B(y,t))} \\
 				& \lesssim \sup_{(y,t) \in \Gamma(x)} \sum_{j=0}^{\infty} 2^{-j(\frac{n}{\tilde{p}}+\eps)} V(y,t)^{-1/\tilde{p}} \norm{f}_{L^{\tilde{p}}(S_j(B(y,t)))}.
 \end{align*}
By application of the doubling condition \eqref{doublingProperty2}, we further get that the above is bounded by a constant times
 \begin{align*}
 				  \sup_{t>0} \sup_{y \in B(x,t)} \sum_{j=0}^{\infty} 2^{-j(\frac{n}{\tilde{p}}+\eps)} \, 2^{j\frac{n}{\tilde{p}}} V(y,2^jt)^{-1/\tilde{p}}  \norm{f}_{L^{\tilde{p}}(B(y,2^jt))} 
 				  \lesssim  \left[\calM(\abs{f}^{\tilde{p}})(x)\right]^{1/\tilde{p}} =  \calM_{\tilde{p}}f(x).
  \end{align*}
 As $\calM_{\tilde{p}}$ is bounded on $L^p(X)$ for every $p \in (\tilde{p},\infty]$, the proof is finished.\\
(ii) First recall that due to Lemma \ref{psi-remark} the operator $e^{-tL}$ can be defined via duality as an operator acting from $L^{\infty}(X)$ to $L^2_{\loc}(X)$ for every $t>0$. With the same reasoning, one can also define for every $p \in (2,\infty)$ via duality $e^{-tL}$ as an operator acting from $L^p(X)$ to $L^2_{\loc}(X)$.\\
Let $p \in (2,\infty]$ and let $f \in L^p(X)$. Then, repeating the arguments in (i), but with the $L^{\tilde{p}}-L^2$ off-diagonal estimates replaced by the Davies-Gaffney estimates for the semigroup, we obtain for every $x \in X$
\begin{align*}
		\calN_{h,L} f(x) 
				& \leq \sup_{(y,t) \in \Gamma(x)} \sum_{j=0}^{\infty} V(y,t)^{-1/2} \norm{ e^{-t^{2m}L} \Eins_{S_j(B(y,t))}f}_{L^2(B(y,t))} \\
				& \lesssim \sup_{(y,t) \in \Gamma(x)} \sum_{j=0}^{\infty} V(y,t)^{-1/2} \exp\left(-\left(\frac{\dist(S_j(B(y,t)),B(y,t))^{2m}}{ct^{2m}}\right)^{\frac{1}{2m-1}}\right)  \norm{f}_{L^2(B(y,2^jt))} \\
				& \lesssim \sup_{t>0} \sup_{y \in B(x,t)} \sum_{j=0}^{\infty} 2^{-j(\frac{n}{2}+\eps)} \, 2^{j\frac{n}{2}} V(y,2^jt)^{-1/2}  \norm{f}_{L^2(B(y,2^jt))}
				 \lesssim \calM_2f(x).
\end{align*}
The claim follows from the fact that $\calM_2$ is bounded on $L^p(X)$ for every $p \in (2,\infty]$.
\end{Proof}

\begin{Remark}
The boundedness of $\calN_{h,L^{\ast}}$ in $L^2(X)$ immediately follows from Lemma \ref{nontang-maxfct} and the assumptions \eqref{H1} and \eqref{L2-Lq-estimate} of \eqref{H3}. 
\end{Remark}

 \begin{Remark} \label{Linfty-average-est}
Let $L$ satisfy \eqref{H1} and \eqref{H2}. Let $p \in (2,\infty]$ and $f \in L^p(X)$.
The proof of Lemma \ref{nontang-maxfct} (ii) in particular shows that for every $t>0$ and every $x \in X$
\begin{align*}
		\abs{A_te^{-t^{2m}L} f(x)}  \leq \frac{1}{V(x,t)} \int_{B(x,t)} \abs{e^{-t^{2m}L} f(y)} \,d\mu(y)
				 \lesssim \calM_2 f(x).
\end{align*} 
The boundedness of $\calM_2$ on $L^p(X)$ for every $p \in (2,\infty]$ then implies that
$
		\norm{A_t e^{-t^{2m}L} f}_{L^p(X)} \lesssim \norm{f}_{L^p(X)}
$
uniformly in $t>0$. 
\end{Remark}

 \begin{Proof}[of Theorem \ref{paraproduct} (i)]
  For $f,g \in L^2(X)$, the Cauchy-Schwarz inequality implies
 \begin{align*}
  	\abs{\skp{\Pi_b(f),g}} 
 		&\leq \left(\iint_{X \times (0,\infty)} \abs{\psi(t^{2m}L) b(x) \cdot A_t(e^{-t^{2m}L}f)(x)}^2 \, \frac{d\mu(x)dt}{t}\right)^{1/2} \\
 		& \qquad \times \left(\iint_{X \times (0,\infty)} \abs{\tilde{\psi}(t^{2m}L^{\ast}) g(x)}^2 \, \frac{d\mu(x)dt}{t}\right)^{1/2}.
 \end{align*}
 The second factor is bounded by a constant times $\norm{g}_{L^2(X)}$ according to assumption \eqref{H1} and \eqref{square-functions}. Recalling the definition of $\nu_{\psi,b}$ in \eqref{CarlesonMeasure}, we see that the first factor is equal to
 \begin{align} \label{parapr-eq1}
     \left(\iint_{X \times (0,\infty)}  \abs{A_t(e^{-t^{2m}L}f)(x)}^2 \, d\nu_{\psi,b}(x,t) \right)^{1/2}.
 \end{align}
 As we assumed $\beta>\frac{n}{4m}$, Proposition \ref{LemmaCarlesonMeasure} yields that $\nu_{\psi,b}$ is a Carleson measure with $\norm{\nu_{\psi,b}}_{\calC}^{1/2} \lesssim \norm{b}_{BMO_L(X)}$.
  On the other hand, observe that the Cauchy-Schwarz inequality yields for every $h \in L^2_{\loc}(X)$ and every $y \in X$ the estimate $\abs{A_th(y)}^2 \leq \frac{1}{V(y,t)} \int_{B(y,t)} \abs{h(z)}^2 \,d\mu(z)$. 
With the help of Theorem \ref{CarlesonDuality}, we can therefore estimate \eqref{parapr-eq1} by a constant times
 \begin{align*}
  	& \norm{\nu_{\psi,b}}_{\calC}^{1/2}  \left(\int_{X} \sup_{(y,t) \in \Gamma(x)} \abs{A_t(e^{-t^{2m}L}f)(y)}^2 \,d\mu(x) \right)^{1/2} \\
  	& \qquad \lesssim \norm{b}_{BMO_L(X)} \left(\int_{X} \sup_{(y,t) \in \Gamma(x)} \frac{1}{V(y,t)} \int_{B(y,t)} \abs{e^{-t^{2m}L}f(z)}^2 \,d\mu(z) \,d\mu(x) \right)^{1/2}\\
  	& \qquad  = \norm{b}_{BMO_L(X)}\norm{\calN_{h,L} f}_{L^2(X)} 
 		\lesssim  \norm{b}_{BMO_L(X)} \norm{f}_{L^2(X)},
 \end{align*}
using the boundedness of $\calN_{h,L}$ on $L^2(X)$ in the last step.
 \end{Proof}

Via the duality of $H^1_{L^{\ast}}(X)$ and $BMO_L(X)$ and with similar arguments as those used in Section 8 of \cite{HofmannMayboroda}, we moreover obtain the following.

\begin{Proof}[of Theorem \ref{paraproduct} (ii), $p=\infty$]
Let $f \in L^{\infty}(X)$. Moreover, let $\eps>0$ and $M \in \N$ with $M>\frac{n}{4m}$ and let $g \in \bH^1_{L^{\ast}}(X)$, where $\bH^1_{L^{\ast}}(X) = H^1_{L^{\ast}}(X) \cap L^2(X)$ as defined in \eqref{Def-bH1}.
For every $R>0$ let us consider $\ell_R$ defined by
\begin{align} \label{Para-BMO-eq0}
	\ell_R(g):=	\skp{\int_{1/R}^R\tilde{\psi}(t^{2m}L) \Eins_{B_R}[\psi(t^{2m}L) b \cdot A_t e^{-t^{2m}L}f]\,\frac{dt}{t},g} ,
\end{align}
where $B_R:=B(0,R)$  and the pairing is that between $H^1_{L^{\ast}}(X)$ and its dual.\\
On the one hand, since $\beta>\frac{n}{4m}$, Theorem \ref{CharacterizationHardy} yields that the function $G$, defined by
\begin{align} \label{Para-BMO-DefG}
		 G(x,t):=\tilde{\psi}(t^{2m}L^{\ast}) g(x), \qquad (x,t) \in X \times (0,\infty),
\end{align}
is an element of $T^1(X)$ with
\begin{align} \label{Para-BMO-eq1}
		\norm{G}_{T^1(X)} = \norm{\scrA G}_{L^1(X)} \lesssim \norm{g}_{H^1_{L^{\ast}}(x)}.
\end{align}
As in the proof before, we use that $\nu_{\psi,b}:=\abs{\psi(t^{2m}L)b(y)}^2 \,\frac{d\mu(y)dt}{t}$ is a Carleson measure with $\norm{\nu_{\psi,b}}^{1/2}_{\calC} \lesssim  \norm{b}_{BMO_L(X)}$. 
Thus, the function $F$, defined by 
\begin{align} \label{Para-BMO-DefF}
		F(x,t):=\psi(t^{2m}L) b(x) \cdot A_t e^{-t^{2m}L} f(x), \qquad (x,t) \in X \times (0,\infty),
\end{align}
is an element of $T^{\infty}(X)$ with 
\begin{align} \label{Para-BMO-eq2} \nonumber
	&\norm{F}_{T^{\infty}(X)}
			 =	\norm{\scrC F}_{L^{\infty}(X)} \\ \nonumber
			& \quad = \norm{ x \mapsto \sup_{B: x \in B} \left(\frac{1}{V(B)} \int_0^{r_B} \int_B \abs{\psi(t^{2m}L)b(y)}^2 \abs{A_te^{-t^{2m}L}f(y)}^2 \,\frac{d\mu(y)dt}{t} \right)^{1/2} }_{L^{\infty}(X)} \\
  			& \quad \lesssim \norm{f}_{L^{\infty}(X)} \norm{\nu_{\psi,b}}^{1/2}_{\calC} 
  					\lesssim \norm{f}_{L^{\infty}(X)} \norm{b}_{BMO_L(X)},
\end{align} 
where we used Remark \ref{Linfty-average-est} in the penultimate step.
This estimate also shows that $\ell_R \in L^2(X)$ for every $R>0$, since Minkowski's inequality, the uniform boundedness of $\{\tilde{\psi}(tL)\}_{t>0}$ and the Cauchy-Schwarz inequality yield
\begin{align*}
		\norm{\ell_R}_{L^2(X)} 
				& = \norm{\int_{1/R}^R \tilde{\psi}(t^{2m}L) \Eins_{B_R} F(\,.\,,t) \,\frac{dt}{t}}_{L^2(X)} 
				 \lesssim \int_{1/R}^R \norm{F(\,.\,,t)}_{L^2(B_R)} \,\frac{dt}{t} \\
				& \leq C_R \left(\int_0^R \int_{B_R} \abs{F(x,t)}^2 \,\frac{d\mu(x)dt}{t} \right)^{1/2} 
				 \leq C_R V(B_R)^{1/2} \norm{F}_{T^{\infty}(X)}.
\end{align*}
Therefore, according to Theorem \ref{CMS}, we obtain from \eqref{Para-BMO-eq1} and \eqref{Para-BMO-eq2}
\begin{align*}
  		\abs{\ell_R(g)}
  					&  \leq \int_0^{\infty}
  								\abs{\skp{\psi(t^{2m}L) b \cdot A_t e^{-t^{2m}L}f, \tilde{\psi}(t^{2m}L^{\ast}) g}} \, \frac{dt}{t}
  					 \lesssim \int_{X} \scrC F(x) \, \scrA G(x) \, d\mu(x) \\
  					& \lesssim \norm{F}_{T^{\infty}(X)} \norm{G}_{T^1(X)} 
  					 \lesssim  \norm{f}_{L^{\infty}(X)} \norm{b}_{BMO_L(X)} \norm{g}_{H^1_{L^{\ast}}(x)}.
 \end{align*}
Since $\bH^1_{L^{\ast}}(X)$ is dense in $H^1_{L^{\ast}}(X)$, the above implies that $\ell_R$ defines a continuous linear functional on $H^1_{L^{\ast}}(X)$ which can, due to Theorem \ref{Duality}, be identified as an element of $BMO_L(X)$ for every $R>0$ with
\begin{align} \label{Para-BMO-eq3}
		\sup_{R>0} \norm{\ell_R}_{BMO_L(X)} \lesssim  \norm{f}_{L^{\infty}(X)} \norm{b}_{BMO_L(X)}.
\end{align}
Moreover, in view of the duality of $T^1(X)$ and $T^{\infty}(X)$ stated in Theorem \ref{CMS}, $\ell_R$ converges pointwise on  $\bH^1_{L^{\ast}}(X)$ for $R \to \infty$ with
\begin{align*}
		\ell_R(g) &= \int_{1/R}^R \skp{\Eins_{B_R} F(\,.\,,t),G(\,.\,,t)} \,\frac{dt}{t} \\
							& \to \int_0^{\infty} \skp{ F(\,.\,,t),G(\,.\,,t)} \,\frac{dt}{t} 
							 = \int_0^{\infty} \skp{\psi(t^{2m}L) b \cdot A_t e^{-t^{2m}L}f,\tilde{\psi}(t^{2m}L^{\ast}) g} \,\frac{dt}{t},  \qquad R \to \infty.
\end{align*}
By uniform boundedness  we can define in this sense $\Pi_b(f)$ as an element of $BMO_L(X)$.
The estimate \eqref{Para-BMO-eq3} finally yields the desired norm estimate of the operator $\Pi_b$.
 \end{Proof}

One possibility to show that $\Pi_b$ also extends to a bounded operator from $L^p(X)$ to $H^p_L(X)$ is the use of the interpolation result for Hardy spaces stated in Proposition \ref{Interpolation}.
We will present a more direct approach, that is similar to the above proof and does not require assumption \eqref{H3}. The idea goes back to \cite{HytoenenWeis}.

\begin{Proof}[of Theorem \ref{paraproduct} (ii), $p \in (2,\infty)$]
Let $\frac{1}{p}+\frac{1}{p'}=1$ and let $f \in L^p(X)$ and $g \in \bH^{p'}_{L^{\ast}}(X)$. 
For every $R>0$, let $\ell_R$ be defined as in \eqref{Para-BMO-eq0}, where the pairing is now that between $H^p_L(X)$ and its dual. Further, let $G$ and $F$ be defined as in \eqref{Para-BMO-DefG} and \eqref{Para-BMO-DefF}.
Then, due to Theorem \ref{Charact-Hp} and the assumption $\tilde{\psi} \in \Psi_{\alpha,\beta}(\Sigma_{\sigma}^0)$ with $\beta>\frac{n}{4m}$, we obtain $G \in T^{p'}(X)$ with
\begin{align} \label{Para-p-eq1}
		\norm{G}_{T^{p'}(X)} = \norm{\scrA G}_{L^{p'}(X)} \lesssim \norm{g}_{H^{p'}_{L^{\ast}}(X)}.
\end{align}
Let us now split $F$ into $F=H \cdot F_0$ with $H(\,.\,,t):=\psi(t^{2m}L) b$ and $F_0(\,.\,,t):=A_t e^{-t^{2m}L}f$. On the one hand, Proposition \ref{LemmaCarlesonMeasure} yields, as before, that $H \in T^{\infty}(X)$ with $\norm{H}_{T^{\infty}(X)} = \norm{\nu_{\psi,b}}_{\scrC}^{1/2} \lesssim \norm{b}_{BMO_{L}(X)}$. Observe that on the other hand $F_0^{\ast}=\calN_{h,L}f$, thus we obtain from Lemma \ref{nontang-maxfct} that $F_0^{\ast} \in L^p(X)$ with $\norm{F_0^{\ast}}_{L^p(X)} \lesssim \norm{f}_{L^p(X)}$.
Therefore, Proposition \ref{Carleson-Tent-Cor} implies that $F \in T^p(X)$ with
\begin{align} \label{Para-p-eq2} \nonumber
		\norm{F}_{T^p(X)} 
						 = \norm{\scrC(H \cdot F_0)}_{L^p(X)} \lesssim \norm{H}_{T^{\infty}(X)} \norm{F_0^{\ast}}_{L^p(X)} 
						 \lesssim \norm{b}_{BMO_{L}(X)} \norm{f}_{L^p(X)}.
\end{align}
Hence, we get due to Theorem \ref{CMS}, H\"{o}lder's inequality and the fact that $\norm{\scrA F}_{L^p(X)}  \lesssim \norm{\scrC F}_{L^p(X)} $ according to \cite{CoifmanMeyerStein}, Theorem 3,
\begin{align*}
		\abs{\ell_R(g)}
					& \leq \int_0^{\infty} \abs{\skp{\tilde{\psi}(t^{2m}L^{\ast})g,\psi(t^{2m}L)b \cdot A_t e^{-t^{2m}L}f}} \,\frac{dt}{t} 
					 \lesssim \int_X \scrA(F)(x) \scrA(G)(x) \,d\mu(x) \\
					& \lesssim \norm{\scrC F}_{L^p(X)} \norm{\scrA G}_{L^{p'}(X)} 
					 \lesssim \norm{b}_{BMO_{L}(X)} \norm{f}_{L^p(X)} \norm{g}_{H^{p'}_{L^{\ast}}(X)},
\end{align*}
where the last step is a consequence of \eqref{Para-p-eq1} and \eqref{Para-p-eq2}.
Since $\bH^{p'}_{L^{\ast}}(X)$ is dense in $H^{p'}_{L^{\ast}}(X)$ and $H^p_L(X)$ was defined as the dual space of $H^{p'}_{L^{\ast}}(X)$, we can therefore identify $\ell_R$ with an element of $H^p_L(X)$. 
With the same reasoning as in the above proof and in view of the duality of $T^p(X)$ and $T^{p'}(X)$, we can finally define $\Pi_b(f)$ as an element of $H^p_L(X)$ and $\Pi_b$ as an operator acting from $L^p(X)$ to $H^p_L(X)$ with
\begin{align*}
		\norm{\Pi_b(f)}_{H^p_L(X)} \leq C  \norm{b}_{BMO_{L}(X)} \norm{f}_{L^p(X)}.
\end{align*}
\end{Proof}

\begin{Remark} \label{Para-identity}
Let us for a moment assume that the semigroup satisfies the conservation property
\begin{align*}
		e^{-tL}(1)=1 \qquad \text{in} \ L^2_{\loc}(X)
\end{align*}
for every $t>0$. 
Let $\psi, \tilde{\psi} \in \Psi(\Sigma_{\sigma}^0)$ and let $g \in H^1_{L^{\ast}}(X)$ be a finite linear combination of $(1,2,M',\eps)$-molecules for some $\eps>0$ and $M' \in \N$ such that the assumptions of Lemma \ref{ReproducingFormula} and Theorem \ref{paraproduct} (ii) are satisfied. 
If one chooses $\psi,\tilde{\psi} \in \Psi(\Sigma_{\sigma}^0)$ such that $\int_0^{\infty} \psi(t) \tilde{\psi}(t) \,\frac{dt}{t}=1$, then Thereom \ref{paraproduct} (ii) implies that $\Pi_b(1) \in BMO_L(X)$ with
\begin{align*}
		\skp{\Pi_b(1),g} 
					&= \int_0^{\infty} \skp{\psi(t^{2m}L) b \cdot  A_t e^{-t^{2m}L}1,\tilde{\psi}(t^{2m}L^{\ast}) g} \,\frac{dt}{t} \\
					&= \int_0^{\infty} \skp{\psi(t^{2m}L) b,\tilde{\psi}(t^{2m}L^{\ast}) g} \,\frac{dt}{t} 
					= \skp{b,g}
\end{align*}
due to the reproducing formula of Lemma \ref{ReproducingFormula}. Since $g$ was arbitrarily chosen from a dense subset of $H^1_{L^{\ast}}(X)$, we thus obtain
\begin{align*}
		\Pi_b(1) = b \qquad \text{in} \ BMO_L(X).
\end{align*}

For the adjoint operator $\Pi^{\ast}_b$ we also obtain, at least at a formal level, the equality 
\begin{align*}
		\Pi^{\ast}_b(1) = \int_0^{\infty} e^{-t^{2m}L^{\ast}} A_t^{\ast} [\overline{\psi(t^{2m}L)b} \cdot \tilde{\psi}(t^{2m}L^{\ast})1] \,\frac{dt}{t} =0,
\end{align*} 
whenever $\tilde{\psi}(tL^{\ast})(1)=0$. 
The condition $\tilde{\psi}(tL^{\ast})(1)=0$ in $L^2_{\loc}(X)$ is fulfilled in the case that $e^{-tL^{\ast}}(1)=1$ in $L^2_{\loc}(X)$ and $\tilde{\psi} \in \Psi_{\beta,\alpha}(\Sigma_{\sigma}^0)$ for some $\alpha>0$ and $\beta>\frac{n}{4m}$, see Lemma \ref{psi-remark}.
\end{Remark}

\subsection{Boundedness of paraproducts - via off-diagonal estimates}
\label{sect-para-further}

Throughout the section we will assume that $L$ satisfies \eqref{H1}, \eqref{H2} and also \eqref{H3}. This is done to avoid technicalities, even if assumption \eqref{H3} will not always be necessary.\\
To obtain further boundedness properties of the paraproduct $\Pi$ defined in \eqref{Paraproduct}, we will consider $\Pi$ in this section as a bilinear operator, initially defined on $L^2(X) \times BMO_L(X)$ for $\psi, \tilde{\psi} \in \Psi(\Sigma_{\sigma}^0)$ by
\begin{align} \label{DefParaproduct-bilinear}
		\Pi(f,g):= \int_0^{\infty} \tilde{\psi}(t^{2m}L)[\psi(t^{2m}L)g \cdot A_t e^{-t^{2m}L}f] \,\frac{dt}{t}
\end{align}
for every $f \in L^2(X)$ and $g \in BMO_L(X)$. 
In Section \ref{paraproduct-L2}, we already showed that $\Pi$ extends to a bounded bilinear operator
\begin{align*}
	&	\Pi: L^2(X) \times BMO_L(X) \to L^2(X), \\
	&	\Pi: L^p(X) \times BMO_L(X) \to H^p_L(X), \qquad 2<p<\infty, \\
	& \Pi: L^{\infty}(X) \times BMO_L(X) \to BMO_L(X),
\end{align*}
if the defining functions of the paraproduct, $\psi,\tilde{\psi} \in \Psi(\Sigma_{\sigma}^0)$, have enough decay at $0$ and infinity, respectively.
In addition, we will now show that $\Pi$ extends to a bounded bilinear operator
\begin{align*}
		& \Pi: L^{\infty}(X) \times H^p_L(X) \to L^p(X), \qquad 1 \leq p < 2, \\
		& \Pi: L^{\infty}(X) \times L^2(X) \to L^2(X), \\
		& \Pi: L^{\infty}(X) \times L^p(X) \to H^p_L(X), \qquad 2 < p < \infty.
\end{align*}

We begin with the simplest case, namely the boundedness of $\Pi: L^{\infty}(X) \times L^2(X) \to L^2(X)$. This is an immediate consequence of quadratic estimates and Remark \ref{Linfty-average-est}.

\begin{Lemma} \label{Paraproduct-L2Linf-bound}
Let $\psi,\tilde{\psi} \in \Psi(\Sigma_{\sigma}^0)$. Then the operator $\Pi$ defined in \eqref{DefParaproduct-bilinear} extends to a bounded operator $\Pi: L^{\infty}(X) \times L^2(X) \to L^2(X)$. I.e. there exists a constant $C>0$ such that for every $f \in L^{\infty}(X)$ and every $g \in L^2(X)$ 
\begin{align*}
		\norm{\Pi(f,g)}_{L^2(X)} \leq C \norm{f}_{L^{\infty}(X)} \norm{g}_{L^2(X)}.
\end{align*}
\end{Lemma}

\begin{Proof}
Let $f \in L^{\infty}(X)$ and $g,h \in L^2(X)$. The Cauchy-Schwarz inequality, Remark \ref{Linfty-average-est} and quadratic estimates for $\{\psi(tL)\}_{t>0}$ and $\{\tilde{\psi}(tL)\}_{t>0}$, which hold due to \eqref{square-functions}, then yield
	\begin{align*}
			 \abs{\skp{\Pi(f,g),h}} 
 			& \leq \left(\int_0^{\infty} \norm{\psi(t^{2m}L)g \cdot A_t e^{-t^{2m}L}f}_{L^2(X)}^2 \,\frac{dt}{t}\right)^{1/2}
								\left(\int_0^{\infty} \norm{\tilde{\psi}(t^{2m}L^{\ast})h}_{L^2(X)}^2 \,\frac{dt}{t}\right)^{1/2}\\
					& \qquad \lesssim \norm{f}_{L^{\infty}(X)} \norm{g}_{L^2(X)} \norm{h}_{L^2(X)}.
	\end{align*}
\end{Proof}

Next, we will show that $\Pi$ extends to a bounded operator $\Pi: L^{\infty}(X) \times H^1_L(X) \to L^1(X)$. We therefore first check that the off-diagonal estimates \eqref{H1-eq1} and \eqref{H1-eq2} of Proposition \ref{H1-bddness} are satisfied.

\begin{Lemma} \label{Para-offdiag}
	Let $\alpha_1,\alpha_2,\beta_1,\beta_2>0$ and let $\psi \in \Psi_{\beta_1,\alpha_1}(\Sigma_{\sigma}^0)$ and $\tilde{\psi} \in \Psi_{\alpha_2,\beta_2}(\Sigma_{\sigma}^0)$. Further, let $\delta>0$ and  $\varphi \in H^{\infty}(\Sigma_{\sigma}^0)$ with $\varphi(z) = \calO(\abs{z}^{\delta})$ for $\abs{z} \to 0$.\\ 
Then for every $\gamma>0$ with $\gamma\leq\min(\beta_1,\alpha_2)$ and $\gamma<\min(\beta_2,\delta)$ there exists some constant $C>0$ such that for every $f \in L^{\infty}(X)$, every $t>0$, arbitrary open sets $E,F \in X$ and every $g \in L^2(X)$ supported in $E$ 
\begin{align*}
 	\norm{\varphi(t^{2m}L)\Pi(f,g)}_{L^2(F)} &\leq C \left(1+\frac{\dist(E,F)^{2m}}{t^{2m}}\right)^{-\gamma} \norm{f}_{L^{\infty}(X)} \norm{g}_{L^2(E)}.
 \end{align*} 
\end{Lemma}

 \begin{Proof}
According to Lemma \ref{Paraproduct-L2Linf-bound}, we can without restriction assume $\dist(E,F)>t$.
Let us abbreviate $\rho:=\dist(E,F)$.
Similar to the proof of \cite{HofmannMartell}, Lemma 2.3, we  define $G_1:=\{x \in X \,:\, \dist(x,F)<\frac{\rho}{2}\}$ and $G_2:=\{x \in X \,:\, \dist(x,F)<\frac{\rho}{4}\}$ and then split $X$ into $X=\bar{G}_2 \cup X \setminus \bar{G}_2$. By construction $G_1,G_2$ are open with $\dist(E,G_1) \geq \frac{\rho}{2}$ and $\dist(F,X \setminus \bar{G}_2) \geq \frac{\rho}{4}$. 
We then obtain via Minkowski's inequality
\begin{align*}
	 \norm{\varphi(t^{2m}L)\Pi(f,g)}_{L^2(F)} 
			&  \leq \int_0^{\infty} \norm{\varphi(t^{2m}L) \tilde{\psi}(s^{2m}L) \Eins_{\bar{G}_2}[\psi(s^{2m}L)g \cdot A_s e^{-s^{2m}L} f]}_{L^2(F)} \,\frac{ds}{s} \\
			& \qquad + \int_0^{\infty} \norm{\varphi(t^{2m}L) \tilde{\psi}(s^{2m}L) \Eins_{X \setminus \bar{G}_2} [\psi(s^{2m}L)g \cdot A_s e^{-s^{2m}L} f]}_{L^2(F)} \,\frac{ds}{s} \\
			&  =: J_{\bar{G}_2} + J_{X \setminus \bar{G}_2}.
\end{align*}
To handle $J_{X \setminus \bar{G}_2}$, we moreover split the integral into two parts $J^1_{X \setminus \bar{G}_2}$ and $J^2_{X \setminus \bar{G}_2}$, representing the integration over $(0,t)$ and $(t,\infty)$, respectively. \\
Observe that due to Proposition \ref{H-inf-offdiag} the operator family $\{\varphi(tL)\tilde{\psi}(sL)\}_{s,t>0}$ satisfies off-diagonal estimates in $s$ of order $\alpha_2$. Using in addition the uniform boundedness of $\{\psi(sL)\}_{s>0}$  on $L^2(X)$ and of $\{A_s e^{-s^{2m}L}\}_{s>0}$  on $L^{\infty}(X)$ in the second step and the substitution $u=\frac{s}{t}$ in the third step, we can therefore estimate the term $J^1_{X \setminus \bar{G}_2}$ by
\begin{align} \label{Para-eq2} \nonumber
 	J^1_{X \setminus \bar{G}_2}
		& \lesssim \int_0^t \left(1+\frac{\dist(F,X \setminus \bar{G}_2)^{2m}}{s^{2m}}\right)^{-\alpha_2} 
				\norm{\psi(s^{2m}L)g \cdot A_s e^{-s^{2m}L} f}_{L^2(X \setminus \bar{G}_2)} \,\frac{ds}{s} \\ \nonumber
		& \lesssim \left(\frac{\dist(E,F)^{2m}}{t^{2m}}\right)^{-\alpha_2} \int_0^t \left(\frac{s}{t}\right)^{2m\alpha_2} \,\frac{ds}{s} \norm{f}_{L^{\infty}(X)} \norm{g}_{L^2(E)} \\
		& \lesssim \left(1+\frac{\dist(E,F)^{2m}}{t^{2m}}\right)^{-\alpha_2}  \norm{f}_{L^{\infty}(X)} \norm{g}_{L^2(E)}.
\end{align}
For an estimate of the second part $J^2_{X \setminus \bar{G}_2}$, let us write for $a>0$
\begin{align} \label{Para-offdiag-eq3}
 	\varphi(tL)\tilde{\psi}(sL) = \left(\frac{t}{s}\right)^a (tL)^{-a} \varphi(tL) (sL)^a \tilde{\psi}(sL).
\end{align}
By assumption on $\varphi$ and $\tilde{\psi}$ there holds $z \mapsto z^{-a}\varphi(z) \in H^{\infty}(\Sigma_{\sigma}^0)$ and $z \mapsto z^a\tilde{\psi}(z) \in \Psi_{\alpha_2+a,\beta_2-a}(\Sigma_{\sigma}^0)$ for every $a>0$ with $a \leq \delta$ and $a<\beta_2$. The application of Proposition \ref{H-inf-offdiag} therefore yields that the operator family $\{(tL)^{-a} \varphi(tL) (sL)^a \tilde{\psi}(sL)\}_{s,t>0}$ satisfies off-diagonal estimates in $s$ of order $\alpha_2+a$ (thus, in particular of order $\alpha_2$). 
Hence, with similar arguments as before, we get
\begin{align} \label{Para-eq4} \nonumber
 	J^2_{X \setminus \bar{G}_2}
			& \lesssim \int_t^{\infty} \left(\frac{t}{s}\right)^{2ma} \left(1+\frac{\dist(F,X \setminus \bar{G}_2)^{2m}}{s^{2m}}\right)^{-\alpha_2} \norm{\psi(s^{2m}L)g \cdot A_s e^{-s^{2m}L} f}_{L^2(X \setminus \bar{G}_2)} \,\frac{ds}{s} \\ 
			& \lesssim \int_t^{\infty} \left(\frac{t}{s}\right)^{2ma} \left(1+\frac{\dist(E,F)^{2m}}{s^{2m}}\right)^{-\alpha_2}  \,\frac{ds}{s}  \norm{f}_{L^{\infty}(X)} \norm{g}_{L^2(E)}. 
\end{align}
Recall that we assumed $\gamma<\min(\beta_2,\delta)$. Thus, we can fix some $a>\gamma$ with $a \leq \delta$ and $a<\beta_2$. For such a choice of $a$ we further get in view of the assumptions $\dist(E,F)>t$ and $\gamma \leq \alpha_2$ 
\begin{align} \label{Para-eq5} \nonumber
 	& \int_t^{\infty} \left(\frac{t}{s}\right)^{2ma} \left(1+\frac{\dist(E,F)^{2m}}{s^{2m}}\right)^{-\alpha_2}  \,\frac{ds}{s} 
		 \leq  \left(\frac{\dist(E,F)^{2m}}{t^{2m}}\right)^{-\gamma}
						\int_t^{\infty} \left(\frac{t}{s}\right)^{2m(a-\gamma)} \,\frac{ds}{s} \\
		& \qquad = \left(\frac{\dist(E,F)^{2m}}{t^{2m}}\right)^{-\gamma} \int_1^{\infty} u^{-2m(a-\gamma)} \frac{du}{u} 
		\lesssim \left(1+\frac{\dist(E,F)^{2m}}{t^{2m}}\right)^{-\gamma}.		
\end{align} 
Combining the equations \eqref{Para-eq2}, \eqref{Para-eq4} and \eqref{Para-eq5} yields the desired estimate for $J_{X \setminus \bar{G}_2}$.\\
Let us now turn to $J_{\bar{G}_2}$. By functional calculus, we obtain from \eqref{Para-offdiag-eq3} that there exists a constant $C>0$ such that for all $s,t>0$
\begin{align*}
 	\norm{\varphi(tL)\tilde{\psi}(sL)}_{L^2(X) \to L^2(X)} \leq C \min\left(1,\frac{t}{s}\right)^{a}.
\end{align*}
Due to the fact that $\bar{G}_2 \subseteq G_1$ and using that $\{\psi(sL)\}_{s>0}$ satisfies off-diagonal estimates in $s$ of order $\beta_1$ according to Proposition \ref{H-inf-offdiag}, we thus obtain
\begin{align} \label{Para-eq6} \nonumber
 	J_{\bar{G}_2} & \lesssim \int_0^{\infty} \min\left(1,\frac{t}{s}\right)^{2ma} \norm{\psi(s^{2m}L)g \cdot A_s e^{-s^{2m}L} f}_{L^2(G_1)} \,\frac{ds}{s} \\
		& \lesssim \int_0^{\infty} \min\left(1,\frac{t}{s}\right)^{2ma} \left(1+\frac{\dist(E,G_1)^{2m}}{s^{2m}}\right)^{-\beta_1} \,\frac{ds}{s} \norm{f}_{L^{\infty}(X)} \norm{g}_{L^2(E)}
\end{align}
Since we assumed $\gamma \leq \beta_1$ and chose $a>\gamma$, we can further estimate the integral in \eqref{Para-eq6} by
\begin{align} \label{Para-eq7} \nonumber
 	& \int_0^{\infty}\min\left(1,\frac{t}{s}\right)^{2ma} \left(1+\frac{\dist(E,F)^{2m}}{s^{2m}}\right)^{-\beta_1} \,\frac{ds}{s} \\ \nonumber
		& \qquad \leq \int_0^{\infty} \min\left(1,\frac{t}{s}\right)^{2ma} \left(\frac{t}{s}\right)^{-2m\gamma} \left(\frac{\dist(E,F)^{2m}}{t^{2m}}\right)^{-\gamma} \,\frac{ds}{s} \\ \nonumber
		& \qquad = \left(\frac{\dist(E,F)^{2m}}{t^{2m}}\right)^{-\gamma} \left[ \int_0^t \left(\frac{s}{t}\right)^{2m\gamma} \,\frac{ds}{s} + \int_t^{\infty} \left(\frac{t}{s}\right)^{2m(a-\gamma)} \,\frac{ds}{s} \right] \\
		& \qquad \lesssim \left(1+\frac{\dist(E,F)^{2m}}{t^{2m}}\right)^{-\gamma}.
\end{align}
The combination of \eqref{Para-eq6} and \eqref{Para-eq7} then gives the desired estimate for $J_{\bar{G}_2}$.
 \end{Proof}

By application of Proposition \ref{H1-bddness} and via interpolation and duality we obtain the following.

\begin{Theorem} \label{para-bdd-hp}
Let $\alpha_1>0$ and $\alpha_2,\beta_1,\beta_2>\frac{n}{4m}$. \\
(i) Let $p \in [1,2)$. If  $\psi \in \Psi_{\beta_1,\alpha_1}(\Sigma_{\sigma}^0)$ and $\tilde{\psi} \in \Psi_{\alpha_2,\beta_2}(\Sigma_{\sigma}^0)$, then the operator $\Pi$ defined in \eqref{DefParaproduct-bilinear} extends to a bounded operator $\Pi: L^{\infty}(X) \times H^p_L(X) \to L^p(X)$. I.e. there exists a constant $C>0$ such that for every $f \in L^{\infty}(X)$ and every $g \in H^p_L(X)$
\begin{align*}
		\norm{\Pi(f,g)}_{L^p(X)} \leq C \norm{f}_{L^{\infty}(X)} \norm{g}_{H^p_L(X)}.
\end{align*}
(ii) Let $p \in (2,\infty)$. If  $\psi \in \Psi_{\alpha_2,\beta_2}(\Sigma_{\sigma}^0)$ and $\tilde{\psi} \in \Psi_{\beta_1,\alpha_1}(\Sigma_{\sigma}^0)$, then the operator $\Pi$ defined in \eqref{DefParaproduct-bilinear} extends to a bounded operator $\Pi: L^{\infty}(X) \times L^p(X) \to H^p_L(X)$. I.e. there exists a constant $C>0$ such that for every $f \in L^{\infty}(X)$ and every $g \in L^p(X)$
\begin{align*}
		\norm{\Pi(f,g)}_{H^p_L(X)} \leq C \norm{f}_{L^{\infty}(X)} \norm{g}_{L^p(X)}.
\end{align*}
\end{Theorem}

\begin{Proof}
Concerning (i), observe that Lemma \ref{Para-offdiag} yields the required off-diagonal estimates for Proposition \ref{H1-bddness}. To see this, choose some $M \in \N$ with $M>\frac{n}{4m}$ and define $\varphi \in H^{\infty}(\Sigma_{\sigma}^0)$ by either $\varphi(z)=(1-e^{-z})^M$ or $\varphi(z)=(ze^{-z})^M$. In both cases, $\abs{\varphi(z)} \lesssim \abs{z}^M$ for  $z \in \Sigma_{\sigma}^0$ with $\abs{z} \leq 1$. Thus, we can choose some $\gamma>\frac{n}{4m}$ with $\gamma \leq \min(\beta_1,\alpha_2)$ and $\gamma<\min(\beta_2,M)$. Due to Lemma \ref{Para-offdiag} the operator family $\{\varphi(t^{2m}L) \Pi(f,g)\}_{t>0}$ satisfies $L^2$ off-diagonal estimates of order $\gamma$ with constant $C \norm{f}_{L^{\infty}(X)}$ for some $C>0$ independent of $f$. We therefore obtain from Proposition \ref{H1-bddness} that $\Pi(f,\,.\,)$ extends to a bounded operator from $H^1_L(X)$ to $L^1(X)$ with 
\begin{align*}
		\norm{\Pi(f,g)}_{L^1(X)} \leq C \norm{f}_{L^{\infty}(X)} \norm{g}_{H^1_L(X)},
\end{align*}
for all $g \in H^1_L(X)$ and some constant $C>0$ independent of $f$ and $g$.
Hence, $\Pi$ extends to a bounded operator $\Pi: L^{\infty}(X) \times H^1_L(X) \to L^1(X)$. Via complex interpolation between $H^1_L(X)$ and $H^2_L(X)=L^2(X)$, which holds due to Proposition \ref{Interpolation}, and interpolation between $L^1(X)$ and $L^2(X)$, we also obtain that $\Pi$ extends to a bounded operator $\Pi: L^{\infty}(X) \times H^p_L(X) \to L^p(X)$ for every $p \in (1,2)$.\\
The assertion (ii) is now obtained from (i) via duality. If $p'$ denotes the conjugate exponent of $p \in (2,\infty)$, then $H^p_L(X)$ was defined as the dual space of $H^{p'}_{L^{\ast}}(X)$. Observe that the dual operator of $\Pi(f,\,.\,)$ is the operator
\begin{align*}
				h \mapsto \int_0^{\infty} \psi(t^{2m}L^{\ast})[\tilde{\psi}(t^{2m}L^{\ast})h \cdot \overline{A_te^{-t^{2m}L} f}] \,\frac{dt}{t},
\end{align*}
which is according to (i) bounded from $H^{p'}_{L^{\ast}}(X)$ to $L^{p'}(X)$ with its operator norm bounded by a constant times $\norm{f}_{L^{\infty}(X)}$. Thus, $\Pi(f,\,.\,)$ is bounded from $L^p(X)$ to $H^p_L(X)$ with
\begin{align*}
		\norm{\Pi(f,g)}_{H^p_L(X)} \leq C \norm{f}_{L^{\infty}(X)} \norm{g}_{L^p(X)}.
\end{align*}
\end{Proof}


\subsection{Leibniz-type rules}

Let us conclude the section with an observation on differentiability properties of paraproducts constructed via functional calculus.
One of the fundamental properties of paraproducts, as they were e.g. considered in \cite{Bony} and \cite{CoifmanMeyer2} in the context of paradifferential operators, is that they satisfy a Leibniz-type rule and ``preserve'' Sobolev classes. 
We will show a corresponding result for the paraproduct $\Pi$ defined in Section \ref{sect-para-further}, according to the general philosophy, ``differentiability'' is not measured in terms of derivatives, but in terms of fractional powers of the operator $L$.\\

Let $\psi, \tilde{\psi} \in \Psi(\Sigma_{\sigma}^0)$. 
Let us recall the paraproduct operator $\Pi$, now more precisely denoted by $\Pi_{\tilde{\psi},\psi}$,  as defined in \eqref{DefParaproduct-bilinear}: For $f \in L^{\infty}(X)$ and $g \in L^2(X)$ we set
\begin{align*}
		\Pi_{\tilde{\psi},\psi}(f,g):= \int_0^{\infty} \tilde{\psi}(t^{2m}L)[\psi(t^{2m}L)g \cdot A_t e^{-t^{2m}L}f] \,\frac{dt}{t}.
\end{align*}

Then the following fractional Leibniz-type rule for paraproducts is valid. 

\begin{Prop} \label{Diff-Para}
	Let $s>0$, let  $\tilde{\psi} \in \Psi_{\beta,\alpha}(\Sigma_{\sigma}^0)$ and $\psi \in \Psi_{\alpha,\beta}(\Sigma_{\sigma}^0)$ for some $\alpha>\frac{s}{2m}$ and $\beta>0$. For $f \in L^{\infty}(X)$ and $g \in \calD(L^{s/2m})$
\begin{align*}
			L^{s/2m} \Pi_{\tilde{\psi},\psi}(f,g) =  \Pi_{\tilde{\psi}_{s},\psi_{s}}(f,L^{s/2m}g),
\end{align*}
where $\tilde{\psi}_{s},\psi_{s}$ are defined by  $\tilde{\psi}_{s}(z):=z^{s/2m}\tilde{\psi}(z)$ and $\psi_{s}(z):=z^{-s/2m}\psi(z)$. \\
Moreover, there exists some constant $C>0$ such that for all $f \in L^{\infty}(X)$ and all $g \in \calD(L^{s/2m})$
\begin{align*}
		\norm{L^{s/2m} \Pi(f,g)}_{L^2(X)}
		\lesssim  \norm{f}_{L^{\infty}(X)}\norm{L^{s/2m}g}_{L^2(X)}.
\end{align*}
\end{Prop}

\begin{Proof}
Due to functional calculus, the proposition is a consequence of the simple calculation
\begin{align*}
		L^{s/2m} \Pi_{\tilde{\psi},\psi}(f,g) & = \int_0^{\infty} (t^{2m}L)^{s/2m} \tilde{\psi}(t^{2m}L)[(t^{2m}L)^{-s/2m} \psi(t^{2m}L)L^{s/2m} g \cdot A_t e^{-t^{2m}L}f] \,\frac{dt}{t} \\
			& = \Pi_{\tilde{\psi}_{s},\psi_{s}}(f,L^{s/2m}g),
\end{align*}
combined with Lemma \ref{Paraproduct-L2Linf-bound}.
\end{Proof}

In view of Theorem \ref{para-bdd-hp}, one can obviously obtain a similar result for the spaces $H^p_L(X)$ and $L^p(X)$, where $p \neq 2$. We refer the reader to Section 8.4 of \cite{HofmannMayborodaMcIntosh} for a discussion of Hardy-Sobolev spaces associated to a second order elliptic operator $L$ in divergence form.\\
A corresponding result for paraproducts constructed via convolution operators is stated in \cite{Christ}, Proposition III.23.\\

With the help of paraproducts and under some additional assumptions on $L$, one can also show a fractional Leibniz-type rule for products of functions. 
It can be understood as a generalization of an inequality of Kato and Ponce, see \cite{KatoPonce}, Lemma X4, where fractional derivatives are replaced by fractional powers of the operator $L$.\\
To simplify notation, we only cite the result for the case $X=\R^n$. For the same result in more general spaces of homogeneous type and a proof of the result, we refer the reader to \cite{FreyHytoenen}.
The essential idea in the proof is a representation the product of two functions with the help of paraproducts. 
That is, via functional calculus one can write 
\begin{align} \label{product-decomp}
 	f \cdot g = \Pi_1(f,g) + \Pi_2(f,g) + \Pi_2(g,f),
\end{align}
where $\Pi_1$ and $\Pi_2$ are appropriately defined paraproduct operators.

\begin{Theorem} \label{LeibnizThm}
Let $L$ satsify (H1) and (H2) and let $e^{-tL}: L^{\infty}(\R^n) \to L^{\infty}(\R^n)$ be bounded uniformly in $t>0$. 
Additionally, let $e^{-tL}(1)=1$  and assume that $\nabla L^{-1/2m}: L^2(\R^n) \to L^2(\R^n)$ is bounded.
Then for every $s \in (0,1)$ there exists some $C>0$ such that for all $f,g \in \calD(L^{s/2m}) \cap L^{\infty}(X)$
\[
 	\norm{L^{s/2m}(fg)}_{L^2(\R^n)} \leq C \norm{L^{s/2m}f}_{L^2(\R^n)} \norm{g}_{L^{\infty}(\R^n)} + C \norm{f}_{L^{\infty}(\R^n)} \norm{L^{s/2m}g}_{L^2(\R^n)}.
\]
\end{Theorem}

\addcontentsline{toc}{section}{References}
\small{

}

\small{\textsc{Institute for Analysis, Karlsruhe Institute of Technology (KIT), Kaiserstr. 89, D-76128 Karlsruhe, Germany} \\ \textit{E-mail address:} \texttt{dorothee.frey@kit.edu}}

\end{document}